%% file: main.tex
\documentclass[11pt]{cocv}
\usepackage{amsmath,amssymb,pstricks,color}
\usepackage{graphicx,epsfig}
\usepackage{coverage}
\usepackage{amsthm}
\usepackage{psfrag}

\theoremstyle{definition}
\newtheorem{rmrks}[thrm]{Remarks}

\newcommand{\real}{{\mathbb{R}}}

\newcommand{\sphere}{{\mathbb{S}}}
\renewcommand{\natural}{{\mathbb{N}}}
\newcommand{\eps}{\epsilon}

\newcommand{\pder}[2]{\frac{\partial #1}{\partial #2}}
\newcommand{\der}[2]{\frac{d #1}{d #2}}

\newcommand{\until}[1]{\{1,\dots, #1\}}

\newcommand{\subscr}[2]{#1_{\textup{#2}}}
\newcommand{\setdef}[2]{\left\{#1 \; | \; #2\right\}}

\newcommand{\ov}{\overline}
\newcommand{\arc}{\operatorname{arc}}
\newcommand{\area}[1]{\operatorname{area}_{\phi}(#1)}
\newcommand{\indicator}[1]{\operatorname{1}_{#1}}
\newcommand{\diam}[1]{\operatorname{diam}(#1)}
\newcommand{\map}[3]{#1: #2 \rightarrow #3}
\newcommand{\dist}{\operatorname{dist}}
\newcommand{\union}{\operatorname{\cup}}
\newcommand{\intersection}{\ensuremath{\operatorname{\cap}}}
\newcommand{\intersect}{\ensuremath{\operatorname{\cap}}}
\newcommand{\vertexes}{{\cal U}}

\newcommand{\PP}{{\cal P}}
\newcommand{\Sc}{{\cal S}}
\newcommand{\EE}{{\mathcal{E}}}
\newcommand{\eqdef}{\triangleq}
\newcommand{\diag}{\operatorname{diag}}

\newcommand\oprocendsymbol{\hbox{$\bullet$}}
\newcommand\oprocend{\relax\ifmmode\else\unskip\hfill\fi\oprocendsymbol}
\def\eqoprocend{\tag*{$\bullet$}}

\begin{document}
\title{Spatially-distributed 
  coverage optimization and control 
  with limited-range interactions}\thanks{Complete preprint version
  with all figures available at http:/\!/motion.csl.uiuc.edu}
\author{Jorge Cort\'es}\address{Coordinated Science Laboratory, University
  of Illinois at Urbana-Champaign, 1308 West Main Street, Urbana, Illinois
  61801, USA;
  \email{\{jcortes,smartine,bullo\}@uiuc.edu}}%
\author{Sonia Mart{\'\i}nez}\sameaddress{1}%
\author{Francesco Bullo}\sameaddress{1}
\date{Submitted on: January 20, 2004. This version: \today}

\begin{abstract}
  This paper presents coordination algorithms for groups of
  mobile agents performing deployment and coverage tasks. As an
  important modeling constraint, we assume that each mobile agent has
  a limited sensing/communication radius.
  Based on the geometry of Voronoi partitions and proximity graphs, we
  analyze a class of aggregate objective functions and propose coverage
  algorithms in continuous and discrete time.
  These algorithms have convergence guarantees and are spatially
  distributed with respect to appropriate proximity graphs.  Numerical
  simulations illustrate the results.
\end{abstract}
%
%
\maketitle
\section*{Introduction}

The current technological development of relatively inexpensive
communication, computation, and sensing devices has lead to an intense
research activity devoted to the distributed control and coordination of
networked systems.  In robotic settings, the study of large groups of
autonomous vehicles is nowadays a timely concern.  The potential advantages
of networked robotic systems are their versatility and robustness in the
realization of multiple tasks such as manipulation in hazardous
environments, pollution detection, estimation and map-building of partially
known or unknown environments.

A fundamental problem in the distributed coordination of mobile robots is
that of providing stable and decentralized control laws that are scalable
with the number of agents in the network. Indeed, since the initial works
from the robotics and ecology communities on similar problems on swarms and
flocking~\cite{CWR:87,AO:86,RCA:98}, there have been various efforts to
provide rigorous procedures with convergence guarantees using a combination
of potential energy shaping methods, gyroscopic forces, and graph
theory~\cite{AJ-JL-ASM:02,ROS-RMM:03c,HT-AJ-GJP:03c,PO-EF-NEL:03,KMP:04}.
In our previous work~\cite{JC-SM-TK-FB:02j,JC-FB:02m}, we studied
distributed algorithms for deployment and optimal coverage problems using
tools from computational geometry, nonsmooth analysis and geometric
optimization. The great interest in coordination problems can be easily
detected in the proceedings of the most recent IEEE Conference on Decision
and Control, the Conference on Cooperative Control and Optimization, or the
International Conference on Distributed Autonomous Robotic Systems.

In devising useful coordination algorithms it is important to progressively
account for the various restrictions that real-world systems impose.
Building on our previous work~\cite{JC-SM-TK-FB:02j,JC-FB:02m}, this paper
develops spatially-distributed algorithms for coverage control amenable to
implementation on (more) realistic models of networks; we do this by
considering the following new aspects.  Firstly, we enforce the
communication or sensing capacity of an agent to be restricted to a bounded
region, typically much smaller than the region where the entire network is
confined. In other words, we assume that the agents will have limited-range
communication and/or sensing capabilities: we refer to these information
exchanges between agents as ``limited-range interactions.''  Secondly, we
provide gradient ascent control laws in both continuous and discrete-time
settings, and we prove that the induced dynamical systems are convergent.
Discrete-time feedback algorithms are indeed the ones truly amenable to
implementation in a group of agents exchanging information over a
communication network.  To deal with these problems, we use a seemingly
unrelated combination of tools from graph theory~\cite{RD:00}, locational
optimization~\cite{AO-BB-KS-SNC:00,ZD-HWH:01}, and systems
theory~\cite{UH-JBM:94}.

The contributions of the paper are the following:
\begin{enumerate}
\item Based on the notion of proximity graph~\cite{JWJ-GTT:92}, we provide
  a formal notion of spatially-distributed vector fields and functions; we
  introduce a novel proximity graph, called \emph{limited-range Delaunay
    graph}, related to the notion of Delaunay graph and disk graph; we
  study the properties of the limited-range Delaunay graph and we show, in
  a formal way, that it can be computed in a spatially-distributed fashion.
  
\item We analyze the smoothness properties of an important class of
  objective functions, called multi-center functions, common in locational
  optimization, quantization theory, and geometric optimization. Our
  analysis supersedes the results
  in~\cite{QD-VF-MG:99,AO-AS:97,AO-BB-KS-SNC:00,ZD-HWH:01,RMG-DLN:98,YA:91}.
  One important objective of the analysis is to determine the extent in
  which certain multi-center functions are spatially distributed and with
  respect to which proximity graphs.
  
\item We consider the problem of steering the location of a group of robots
  to local maxima of the objective function. To achieve this objective in
  continuous and discrete-time, we design novel spatially-distributed
  control laws for groups of robots. We formally analyze their performance
  and illustrate their behavior in simulations.
\end{enumerate}

To perform the smoothness analysis in (ii) and the stability analysis in
(iii), we prove useful extensions of the Conservation-of-Mass Law from
fluid dynamics and of the discrete-time LaSalle Invariance Principle,
respectively.  These extensions are, to the best of our knowledge, not
present in classical texts on the subject.

It is worth remarking that one fundamental scientific problem in the
study of coordination algorithms is scalability with respect to
communication complexity. In other words, it is important to design
algorithms with communication requirements that scale nicely (e.g.,
linearly) with the number of agents in the network. However, it is
impossible to quantify the communication complexity of any algorithm
without introducing a detailed communication model.  Adopting a
computational geometric approach, this paper classifies the complexity
of coordination algorithms in terms of the proximity graphs with
respect to which the algorithms are spatially distributed.  The
underlying assumption is that low complexity proximity graphs (e.g.,
graphs with a low number of edges) will require limited communication
in a realistic implementation.

Throughout the paper we shall consider purposefully only extremely
simple models for the dynamics of each individual agent. In
particular, we shall assume that the state of each agent is a point in
$\real^2$ and that the dynamical model of each agent is an integrator
(indeed, we shall interchangeably refer to agent as a location or
point).  This feature is a natural consequence of our focus on
network-wide coordination aspects.

The organization of the paper is as follows.  In
Section~\ref{sec:preliminaries} we review various preliminary mathematical
concepts, and we introduce the notion of proximity graph function and of
spatially-distributed map.  In Section~\ref{sec:locational-optimization} we
study the smoothness of the multi-center function and show in what sense
its partial derivative is spatially distributed.  In
Section~\ref{sec:spat-dist-algorithms} we design spatially-distributed
coverage algorithms, first in continuous-time and then in discrete-time.
The numerical outcomes of the algorithms' implementation are reported in
Section~\ref{sec:simulations}.  Finally, we discuss possible avenues of
future research in Section~\ref{sec:conclusions}.

\section{Preliminaries}
\label{sec:preliminaries}
In this section we present a variety of preliminary concepts.  Graph
theory and proximity graphs from computational geometry are basic
notions that will later allow us to introduce the notion of
spatially-distributed vector fields and algorithms.

\subsection{Basic notions in graph theory}

Here we gather some basic facts on graph theory; for a comprehensive
treatment we refer the reader to~\cite{RD:00}. Given a set
$\vertexes$, recall that $2^{\vertexes}$ is the collection of subsets
of $\vertexes$. A \emph{graph} $\GG = (\vertexes,\EE)$ consists of a
\emph{vertex set} $\vertexes$ and an \emph{edge set} $\EE \subseteq
2^{\vertexes \times \vertexes}$. A graph $(\vertexes,\EE)$ is
\emph{undirected} if $(i,j)\in\EE$ implies $(j,i)\in\EE$.  If $(i,j)
\in \EE$, then vertex $j$ is a \emph{neighbor (in $\GG$)} of vertex
$i$.  Let $\map{\NN_{\GG}}{\vertexes}{2^{\vertexes}}$ associate to the
vertex $i$ the set of its neighbors in $\GG$.  A graph $\GG$ is called
\emph{complete} if any two different vertexes in $\vertexes$ are
neighbors, i.e., $\EE = \vertexes \times \vertexes \setminus \diag
(\vertexes \times \vertexes)$.  This is usually denoted by $K^n$.  A
\emph{path} connecting vertex $i$ to vertex $j$ is a sequence of
vertexes $\{i_0=i,i_1,\dots,i_k,i_{k+1}=j\}$ with the property that
$(i_l,i_{l+1})\in\EE$ for all $l\in\{0,\dots,k\}$. A graph $\GG$ is
\emph{connected} if there exists a path connecting any two vertexes
of~$\GG$.  Given two graphs $\GG_1 = (\vertexes_1,\EE_1)$ and $\GG_2 =
(\vertexes_2,\EE_2)$, the \emph{intersection graph} $\GG_1
\intersection \GG_2$ is the graph $(\vertexes_1 \intersection
\vertexes_2, \EE_1 \intersection \EE_2)$, and the \emph{union graph}
$\GG_1 \union \GG_2$ is the graph $(\vertexes_1 \union \vertexes_2,
\EE_1 \union \EE_2)$.

A graph $\GG_1=(\vertexes_1,\EE_1)$ is a \emph{subgraph} of a graph
$\GG_2 = (\vertexes_2,\EE_2)$ if $\vertexes_1 \subseteq \vertexes_2$
and $\EE_1 \subseteq \EE_2$.  Alternatively, $\GG_2$ is said to be a
\emph{supergraph} of $\GG_1$.  Formally, we set $\GG_1 \subseteq
\GG_2$.  If $\GG_1 \subseteq \GG_2$ and $\GG_1$ contains all the edges
$(i,j) \in \EE_2$ with $i,j \in \vertexes_1$, then $\GG_1$ is called
an \emph{induced subgraph of $\GG_2$}.  A subgraph $\GG_1$ of $\GG_2$
is called \emph{spanning} if $\vertexes_1 = \vertexes_2$.  A
\emph{cycle} of $\GG$ is a subgraph where every vertex has exactly two
neighbors.  An \emph{acyclic} graph is a graph that contains no
cycles. A \emph{tree} is a connected acyclic graph. Given a connected
graph $\GG$, assign to each edge an specific length or \emph{weight}.
The weight of a subgraph of $\GG$ is the sum of the weights of its
edges.  A \emph{minimum spanning tree of $\GG$} is a spanning tree
with the smallest possible weight.  In general, there might exist more
than one minimum spanning tree of $\GG$, all with the same weight.

\subsection{Voronoi partitions and proximity graphs}
\label{sec:geographs}
We start by reviewing the notion of Voronoi partition generated by
sets of points on the Euclidean plane; we refer the reader
to~\cite{MdB-MvK-MO:97,AO-BB-KS-SNC:00} for comprehensive treatments.
Next, we shall present some relevant concepts on proximity graph
functions, that is, on graphs whose vertex set is (in 1-1
correspondence with) a set of distinct points on the plane and whose
edge set is a function of the relative locations of the point set.
This notion is an extension of the notion of proximity graph as
explained in the survey article~\cite{JWJ-GTT:92}; see
also~\cite{XYL:03} and the literature on topology control in wireless
networks for related references.

A \emph{covering} of $\real^2$ is a collection of subsets of $\real^2$
whose union is $\real^2$; a \emph{partition} of $\real^2$ is a
covering whose subsets have disjoint interiors.  Let $\PP$ be a set of
$n$ distinct points $\{p_1,\dots,p_n\}$ in $\real^2$. The
\emph{Voronoi partition of $\real^2$} generated by $\PP$ with respect
to the Euclidean norm is the collection of sets
$\{V_i(\PP)\}_{i\in\until{n}}$ defined by
\begin{equation*}
  V_i(\PP) = \setdef{q\in\real^2}{\|q-p_i\|\leq \|q-p_j\| \, , \;
    \text{for all} \; p_j\in\PP} .
\end{equation*}
Here, $\|\cdot\|$\ denotes the standard Euclidean norm.  It is
customary and convenient to refer to $V_i(\PP)$ as $V_i$. The boundary
of each set $V_i$ is the union of a finite number of segments and
rays.

Let $\Sigma_n$ be the set of permutations of $n$ elements.  A map
$\map{f}{X^n}{2^{X\times X}}$ is \emph{$\Sigma_n$-equivariant} if for
all $(x_1,\dots,x_n)\in X^n$ and $\sigma \in \Sigma_n$, $(x_i,x_j)\in
f(x_1,\dots,x_n)$ implies $(x_{\sigma(i)},x_{\sigma(j)})\in
f(x_{\sigma(1)},\dots,x_{\sigma(n)})$.

A \emph{proximity graph function} associates to a set of $n$ distinct
points $\PP=\{p_1,\dots,p_n \}$ in $\real^2$ a graph with vertex set
$\PP$ and edge set $\EE(p_1,\dots,p_n)$, where
$\map{\EE}{(\real^2)^n}{ 2^{ \real^2\times\real^2 }}$ is a
$\Sigma_n$-equivariant map with the property that
$\EE(p_1,\dots,p_n)\subseteq \PP^2 = \{p_1,\dots,p_n\}^2 =
\{p_1,\dots,p_n\} \times \{p_1,\dots,p_n\}$.

Note that, since the map $\EE$ is $\Sigma_n$-equivariant, the value of
$\EE(p_1,\dots,p_n)$ is independent of the ordering of the elements
$(p_1, \dots, p_n)$, and therefore, with a slight abuse of notation,
we will write it as $\{p_1,\dots,p_n\} =\PP \mapsto \EE(\PP)$, and
refer to it as the \emph{proximity edge function} corresponding to the
proximity graph function $\PP\mapsto\GG(\PP)$.

For $p\in\real^2$ and $r\in\real_+ = [0,+\infty)$, let $B_{r}(p) =
\setdef{q \in \real^2}{\|q-p\|\le r} $ denote the closed ball in
$\real^2$ centered at $p$ of radius $r$.  Now, for $r\in\real_+$, we
have the following proximity graph functions:
\begin{enumerate}
\item the \emph{Delaunay graph} $\PP\mapsto\subscr{\GG}{D}(\PP) =
  (\PP,\subscr{\EE}{D}(\PP))$ has edge set
  \begin{equation*}
  \subscr{\EE}{D}(\PP)=\setdef{(p_i,p_j) \in \PP^2
    \setminus \diag (\PP^2)}{ V_i(\PP)\intersect
    V_j(\PP) \neq  \emptyset} \,;
  \end{equation*}
\item the \emph{$r$-disk graph} $\PP\mapsto\subscr{\GG}{disk}(\PP,r) =
  (\PP,\subscr{\EE}{disk}(\PP,r))$ has edge set
  \begin{equation*}
    \subscr{\EE}{disk}(\PP,r) = \setdef{(p_i,p_j)\in \PP^2 
    \setminus \diag (\PP^2)}{ \| p_i - p_j \| \le r} \,;
  \end{equation*}
\item the \emph{$r$-Delaunay graph}
  $\PP\mapsto\subscr{\GG}{disk\intersection{D}}(\PP,r)$ is the
  intersection of $\subscr{\GG}{disk}(\PP,r)$ and
  $\subscr{\GG}{D}(\PP)$;
  
\item the \emph{$r$-limited Delaunay (or, limited-range Delaunay) graph}
  $\PP\mapsto \subscr{\GG}{LD}(\PP,r) = (\PP,\subscr{\EE}{LD}(\PP,r))$
  consists of the edges $(p_i,p_j) \in \PP^2\setminus \diag (\PP^2)$ with
  the property that
  \begin{align}
    \label{eq:dfn-Voronoi-r}
    \Delta_{ij}(\PP,r) \eqdef \big(V_i(\PP) \cap
    B_{\frac{r}{2}}(p_i)\big) \intersect \big(V_j(\PP) \cap
    B_{\frac{r}{2}}(p_j) \big) \neq \emptyset \,;
  \end{align}
\item the \emph{Gabriel graph}, $\PP \mapsto
  \subscr{\GG}{G}(\PP)=(\PP,\subscr{\EE}{G}(\PP))$ consists of the
  edges $(p_i,p_j) \in \PP^2\setminus \diag (\PP^2)$ with the property
  that
  \begin{equation}\label{eq:dfn-GG}
    p_k \not \in
    \interior B_{\frac{\| p_i - p_j \|}{2}}
    \left(\frac{p_i+p_j}{2}\right) , \quad \text{for all} \;   k \in
    \until{n} \setminus \{i,j\} \, ;
  \end{equation}
\item an \emph{Euclidean Minimum Spanning Tree}, $\PP \mapsto
  \subscr{\GG}{EMST}(\PP) = (\PP,\subscr{\EE}{EMST}(\PP))$ is defined
  as a minimum spanning tree of the complete graph $(\PP, \PP^2
  \setminus \diag (\PP^2)$, whose edge $(p_i,p_j)$ has weight $\| p_i
  - p_j \|$, for $(i,j) \in \until{n}$.
\end{enumerate}

Figure~\ref{fig:example-graphs} presents an example of these proximity
graphs for a random configuration of points.  In general, one can
prove that $ \subscr{\GG}{EMST}(\PP) \subseteq \subscr{\GG}{G}(\PP)
\subseteq \subscr{\GG}{D}(\PP)$ (see for instance~\cite{JWJ-GTT:92}).
While the $r$-Delaunay graph has been studied in earlier
works~\cite{XYL:03,JG-LJG-JH-LZ-AZ:01}, the $r$-limited Delaunay graph
appears not to have been considered.  In the next proposition, we
study some basic useful properties of these graphs.  Before presenting
it, let us recall the following notation from computational geometry.
We denote the cardinality of a set~$S$ by $\# S$. Given $f: \natural
\rightarrow \natural$ and a function $F$ that associates to a set of
$n$ distinct points $\PP=\{p_1,\dots,p_n \}$ in $\real^2$ a
non-negative integer number $F(\{p_1,\dots,p_n\}) \in \natural$, we
denote $F = O (f(n))$ (respectively, $F = \Omega (f(n))$) if and only
if there exists $C \in \ov{\real}_+$ such that $F(\{p_1,\dots,p_n\})
\le C \, f(n)$ (respectively, $F(\{p_1,\dots,p_n\}) \ge C \, f(n)$)
for all distinct $p_1,\dots,p_n \in \real^2$. We denote $F = \Theta
(f(n))$ if and only if both $F = O (f(n))$ and $F = \Omega (f(n))$
hold true.

{
\psfrag{G1}[cc][cc]{{\small Delaunay graph}}
\psfrag{G2}[cc][cc]{{\small $r$-disk graph}}
\psfrag{G3}[cc][cc]{{\small $r$-Delaunay graph}}
\psfrag{G4}[cc][cc]{{\small $r$-limited Delaunay graph}}
\psfrag{G5}[cc][cc]{{\small Gabriel graph}}
\psfrag{G6}[cc][cc]{{\small EMST graph}}

\begin{figure}[htbp]
  \centering
  \fbox{\includegraphics[width=.31\linewidth]{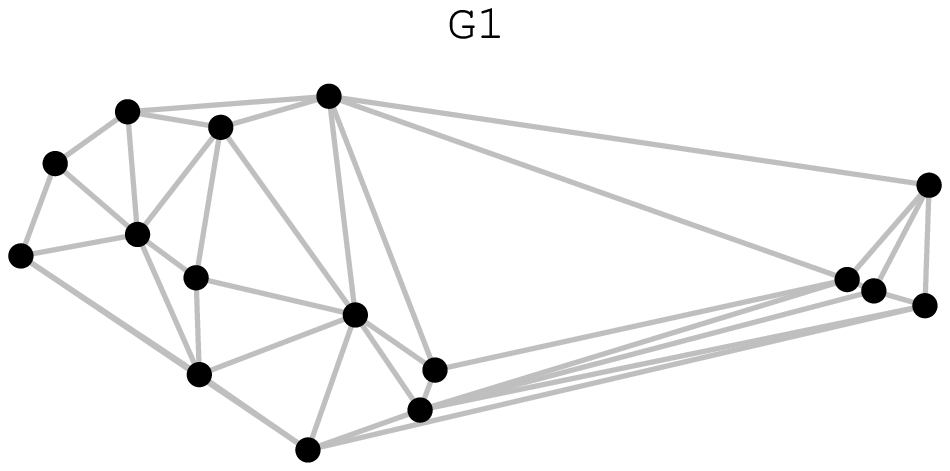}}%
  \fbox{\includegraphics[width=.31\linewidth]{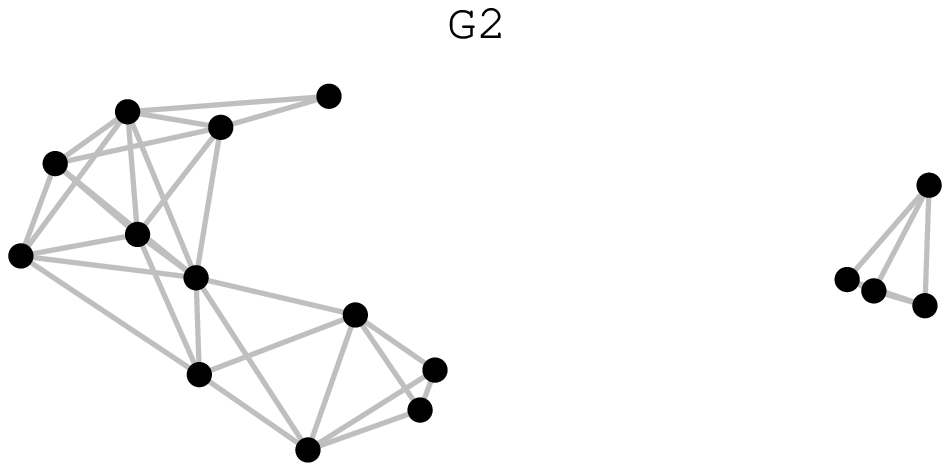}}%
  \fbox{\includegraphics[width=.31\linewidth]{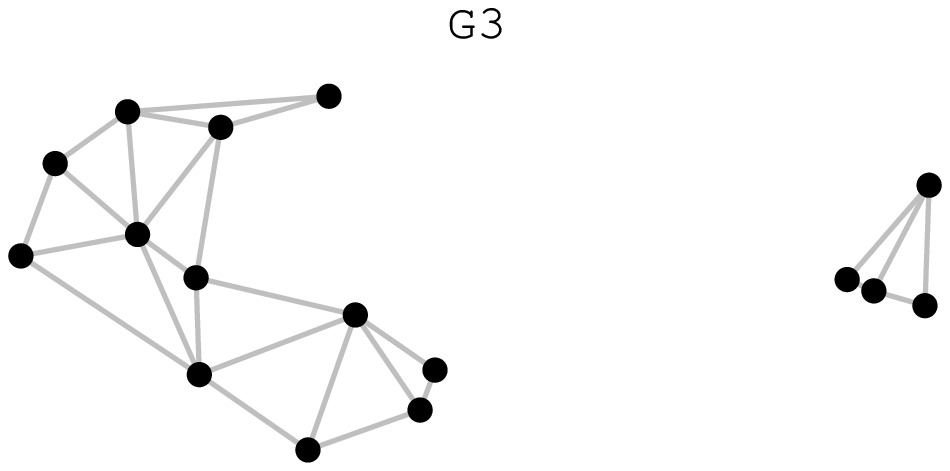}}\\
  \fbox{\includegraphics[width=.31\linewidth]{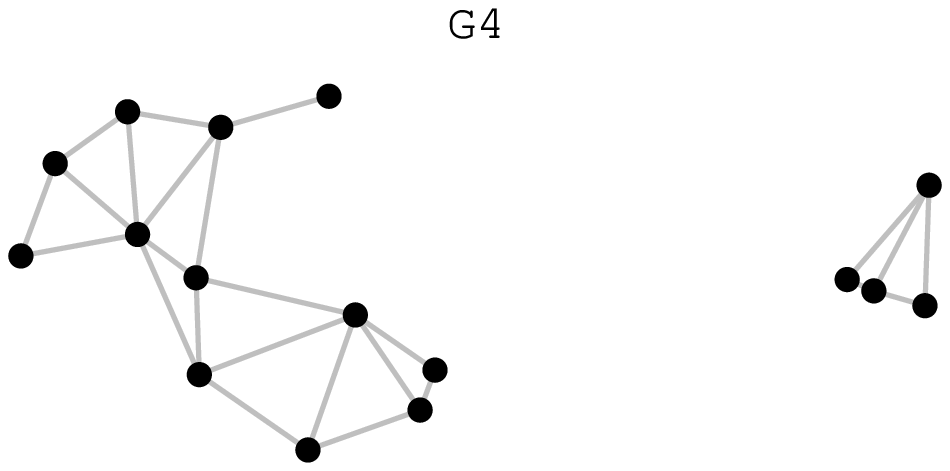}}%
  \fbox{\includegraphics[width=.31\linewidth]{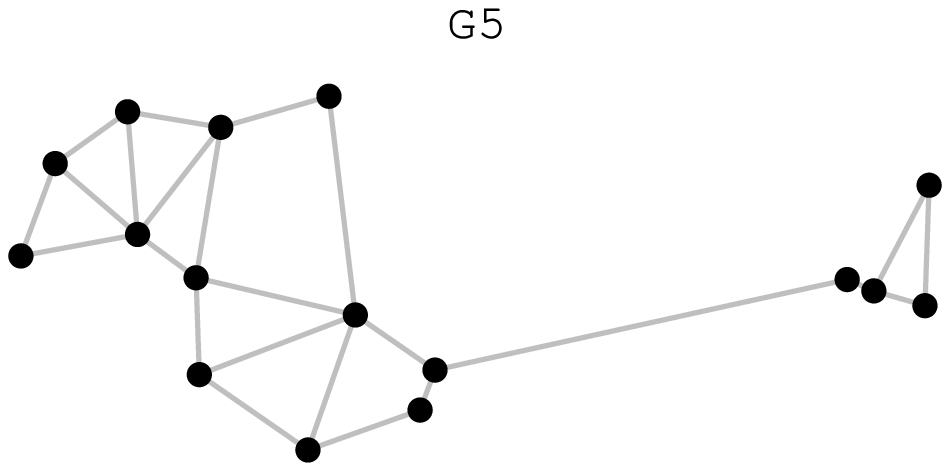}}%
  \fbox{\includegraphics[width=.31\linewidth]{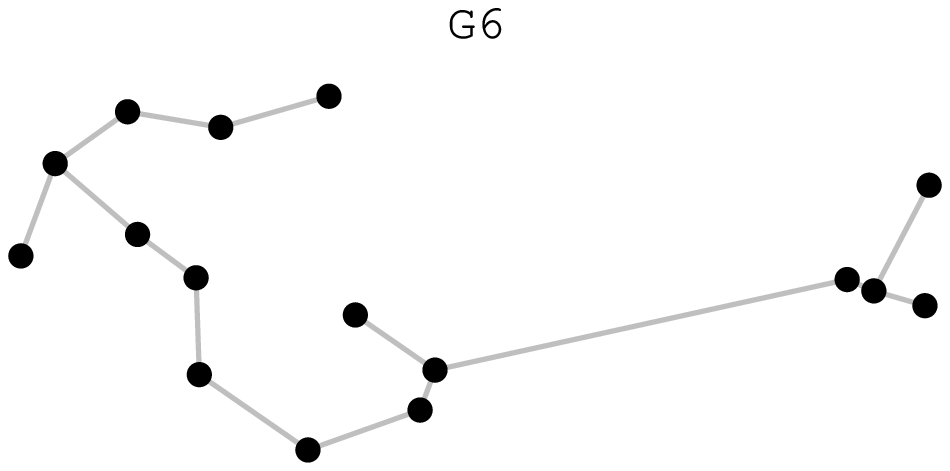}}
  \caption{From left to right, and from up to down, Delaunay,
    $r$-disk, $r$-Delaunay, $r$-limited Delaunay, Gabriel and
    Euclidean Minimum Spanning Tree graphs for a configuration of $16$
    generators with coordinates contained in the rectangle~$[0,1.9]
    \times [0,.75]$. The parameter $r$ is taken equal to
    $.45$.}\label{fig:example-graphs}
\end{figure}
}

\begin{prpstn}\label{prop:graph-theory}
  Let $\PP$ be a set of $n$ distinct points $\{p_1,\dots,p_n\}$ in
  $\real^2$, and let $r \in \real_+$. The following statements hold
  \begin{enumerate}
  \item $\subscr{\GG}{disk\intersection{G}}(\PP,r) \subseteq
    \subscr{\GG}{LD}(\PP,r) \subseteq
    \subscr{\GG}{disk\intersection{D}}(\PP,r)$;
  \item $\subscr{\GG}{disk}(\PP,r)$ is connected if and only if
    $\subscr{\GG}{LD}(\PP,r)$ is connected;
  \item $\# \EE_{\textup{LD}}(\PP) = O(n)$ and $\#
    \EE_{\textup{disk}}(\PP,r) = O(n^2)$. If $\subscr{\GG}{disk}(\PP,r)$
    is connected, then $\# \EE_{\textup{LD}}(\PP) = \Theta(n)$.
  \end{enumerate} 
\end{prpstn}
\begin{proof}
  We first prove the inclusion
  $\subscr{\GG}{disk\intersection{G}}(\PP,r) \subseteq
  \subscr{\GG}{LD}(\PP,r)$.  Let $(p_i,p_j) \in
  \subscr{\EE}{disk\intersection{G}}(\PP,r)$. From the definition of
  the Gabriel graph, we deduce that $\| \frac{p_i+p_j}{2} - p_i \| =
  \| \frac{p_i+p_j}{2} - p_j \| \le \| \frac{p_i+p_j}{2} - p_k \|$,
  for all $k \in \until{n} \setminus \{ i,j\}$, and therefore,
  $\frac{p_i+p_j}{2} \in V_i \cap V_j$. Since $(p_i,p_j) \in
  \subscr{\EE}{disk}(\PP,r)$, we deduce that $\frac{p_i+p_j}{2} \in
  B_{\frac{r}{2}} (p_i) \cap B_{\frac{r}{2}} (p_j)$, and hence
  equation~\eqref{eq:dfn-Voronoi-r} holds, i.e., $(p_i,p_j) \in
  \subscr{\EE}{LD}(\PP,r)$. The second inclusion in (i) is
  straightforward: if $(p_i,p_j) \in \subscr{\EE}{LD}(\PP,r)$, then
  equation~\eqref{eq:dfn-Voronoi-r} implies that $V_i(\PP) \cap V_j
  (\PP) \neq \emptyset$, i.e., $(p_i,p_j) \in \subscr{\EE}{D}(\PP)$.
  Since clearly $(p_i,p_j) \in \subscr{\EE}{disk}(\PP,r)$, we conclude
  (i).  The statement (ii) is a consequence of the following more
  general fact: the $r$-disk graph $\subscr{\GG}{disk}(\PP,r)$ is
  connected if and only if $\subscr{\GG}{EMST}(\PP) \subseteq
  \subscr{\GG}{disk}(\PP,r)$.  The proof of this fact is as follows.
  If $\subscr{\GG}{EMST}(\PP) \subseteq \subscr{\GG}{disk}(\PP,r)$,
  then clearly $\subscr{\GG}{disk}(\PP,r)$ is connected. To prove the
  other implication, assume that $\subscr{\GG}{disk}(\PP,r)$ is
  connected.  We reason by contradiction. Let $\subscr{\GG}{EMST}(\PP)
  \not \subseteq \subscr{\GG}{disk}(\PP,r)$, i.e, there exists $p_i$
  and $p_j$ with $(p_i,p_j) \in \subscr{\EE}{EMST}(\PP)$ and $\| p_i -
  p_j\| > r$.  If we remove this edge from $\subscr{\EE}{EMST}(\PP)$,
  the tree becomes disconnected into two connected components $T_1$
  and $T_2$, with $p_i \in T_1$ and $p_j \in T_2$. Now, since by
  hypothesis the $r$-disk graph $\subscr{\GG}{disk}(\PP,r)$ is
  connected, there must exist $k, l \in \until{n}$ such that $p_k \in
  T_1$, $p_l \in T_2$ and $\|p_k - p_l \| \le r$. If we add the edge
  $(p_k,p_l)$ to the set of edges of $T_1 \cup T_2$, the obtained
  graph $\GG$ is acyclic, connected and contains all the vertexes
  $\PP$, i.e., $\GG$ is a spanning tree.  Moreover, since $\|p_k - p_l
  \| \le r < \| p_i - p_j\|$ and $T_1$ and $T_2$ are induced subgraphs
  of $\subscr{\GG}{EMST}(\PP)$, we conclude that $\GG$ has smaller
  length than $\subscr{\GG}{EMST}(\PP,r)$, which is a contradiction.
  As a consequence, we deduce that if $\subscr{\GG}{disk}(\PP,r)$ is
  connected, then $\subscr{\GG}{EMST}(\PP) \subseteq
  \subscr{\GG}{disk\intersection{G}}(\PP,r)$. Using (i), we conclude
  that $\subscr{\GG}{LD}(\PP,r)$ is connected.  Finally, the statement
  (iii) follows from (i) and~(ii) by noting that $\#
  \EE_{\textup{D}}(\PP) \le 3 n - 6$ (see, for
  instance,~\cite{AO-BB-KS-SNC:00}) and $\# \EE_{\textup{EMST}}(\PP) =
  n-1$.
\end{proof}

Let us make the following observations concerning
Proposition~\ref{prop:graph-theory}.

\begin{rmrks}\label{re:jorge-is-smart}
  As before, let $\PP$ be a set of $n$ distinct points $\{p_1,\dots,p_n\}$
  in $\real^2$, and let $r \in \real_+$.
  \begin{enumerate}
  \item The $r$-Delaunay graph does \emph{not} coincide in general
    with the $r$-limited Delaunay graph.
    Figure~\ref{fig:counter-example} illustrates a point set $\PP$ in
    which $p_l$ is a neighbor of $p_i$ in
    $\subscr{\GG}{disk\intersect{D}}(\PP,r)$ but not in
    $\subscr{\GG}{LD}(\PP,r)$.
    \begin{figure}[htb] 
      \begin{center}
        \resizebox{.32\linewidth}{!}{\input{counter-example.tex}}
      \end{center}
      \caption{Example point set for which the $r$-Delaunay graph strictly
        contains the $r$-limited Delaunay graph: $p_l$ is a neighbor
        of $p_i$ in $\subscr{\GG}{disk\intersect{D}}(\PP,r)$ but not
        in $\subscr{\GG}{LD}(\PP,r)$.}
      \label{fig:counter-example}
    \end{figure}
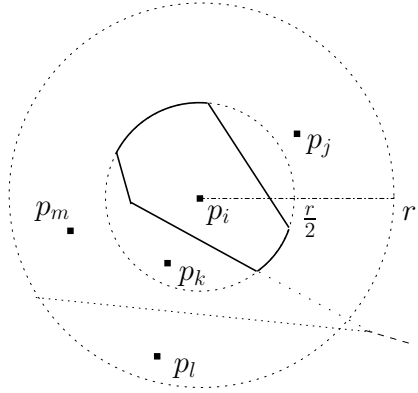              
    
  \item The collection $\{V_i(\PP) \cap
    B_{\frac{r}{2}}(p_i)\}_{i\in\until{n}}$ is a partition of the set
    $\cup_{i} B_{\frac{r}{2}}(p_i) \subset \real^2$.  The boundary of
    each set $V_i(\PP) \cap B_{\frac{r}{2}}(p_i)$, $i\in\until{n}$, is
    the union of a finite number of segments and arcs; see
    Figure~\ref{fig:limited-Delaunay}.  Therefore, at fixed $\PP$,
    there exist $n$ numbers $M_i(r) \ge 0$, $i\in\until{n}$, of
    distinct arcs $\arc_{i,1}(r),\dots,\arc_{i,M_i(r)}(r)$ of radius
    $\frac{r}{2}$ in $\partial (V_i(\PP) \cap B_{\frac{r}{2}}(p_i))$
    with the property that
    \begin{align*}
      \partial \big( V_i(\PP) \cap B_{\frac{r}{2}}(p_i) \big) = 
      \left( \union_{j \in \NN_{{\cal G}_\textup{LD}(\PP,r)}(p_i) }
        \Delta_{ij}(r) \right)
      \union 
      \left( \union_{l \in \until{M_i(r)}} \arc_{i,l}(r) \right),
    \end{align*}
    where we recall that $\NN_{\subscr{\GG}{LD}(\PP,r)}(p_i)$ denotes
    the set of neighbors in $\subscr{\GG}{LD}(\PP,r)$ of the vertex
    $p_i$.
    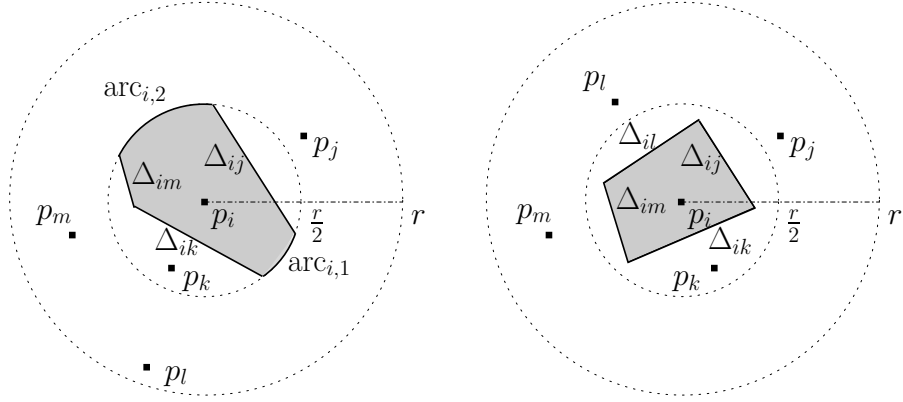
\begin{figure}[htb] 
      \begin{center}
        \resizebox{.32\linewidth}{!}{\input{restricted-Delaunay.tex}}
        \qquad
        \resizebox{.32\linewidth}{!}{\input{restricted-Delaunay-2.tex}}
      \end{center} 
      \caption{The shaded regions are examples of sets  $V_i(\PP) \cap
        B_{\frac{r}{2}}(p_i)$. In the right figure, the boundary of
        the set $V_i(\PP) \cap B_{\frac{r}{2}}(p_i)$ does not contain
        any arc.}
      \label{fig:limited-Delaunay}
    \end{figure}%
  \item \label{item-remark} If $\| p_i-p_j\|$ is strictly greater than
    $r$, then the half plane $\setdef{q\in\real^2}{\|q-p_i\|\leq
      \|q-p_j\|}$\ contains the ball $B_{\frac{r}{2}}(p_i)$.
    Accordingly,
    \begin{align*}
      B_{\frac{r}{2}}(p_i)\cap V_i(\PP) &= B_{\frac{r}{2}}(p_i)\cap
      \setdef{q\in\real^2}{\|q-p_i\|\leq \|q-p_j\| \, , \; \text{for
          all} \;  p_j\in\PP}  \\
      &= B_{\frac{r}{2}}(p_i)\cap \setdef{q\in\real^2}{\|q-p_i\|\leq
        \|q-p_j\| \, , \; \text{for all} \; p_j\in
        \NN_{\subscr{\GG}{disk}(\PP,r)}(p_i) } .
    \end{align*}
  \item It is customary and convenient to refer to the various
    proximity graphs functions without explicitly writing the argument
    $\PP$.  \oprocend
  \end{enumerate}  
\end{rmrks}

Finally, we conclude this section with a general note.

\begin{rmrk}\label{re:Voronoi-covering}
  In the previous definitions, we have emphasized the fact that the points
  $\{p_1,\dots,p_n \}$ are distinct.  Occasionally though, we will consider
  tuples of elements of $\real^2$ of the form $(p_1,\dots,p_n)$, i.e.,
  ordered sets of possibly coincident points.  In this case, it is useful
  to adopt the following notation: given a tuple
  $(p_1,\dots,p_n)\in(\real^2)^n$, possibly containing repeated entries, we
  let $\{p_1,\dots,p_n\}$, or equivalently $\PP$, denote the associated
  point set that only contains the corresponding distinct points.  The
  cardinality of $\PP=\{p_1,\dots,p_n\}$ is therefore less than or equal to
  $n$. More precisely, if $\Sc$ denotes the set
  \begin{equation}
    \label{eq:Sc}
    \Sc = \setdef{(p_1,\dots,p_n) \in (\real^2)^n}{p_i = p_j \;
      \text{for some} \; i, j \in \until{n}, \, i \neq j} \, ,
  \end{equation}
  then $\#\PP < n$ if $(p_1,\dots,p_n) \in \Sc$ and $\#\PP = n$ if
  $(p_1,\dots,p_n) \not \in \Sc$.  The \emph{Voronoi covering $\VV
    (p_1,\dots,p_n) = \{V_i(p_1,\dots,p_n) \}_{i \in \until{n}}$
    generated by the tuple $(p_1,\dots,p_n)$} is defined by assigning
  to each point $p_i$ its corresponding Voronoi cell in the Voronoi
  partition generated by $\PP$. Note that coincident points in the
  tuple $(p_1,\dots,p_n)$ have the same Voronoi cell.  It is
  interesting to note that if the points $p_1,\dots,p_n \in \real^2$
  are not necessarily distinct, then both $\# \EE_{\textup{D}}(\PP) =
  O(n^2)$ and $\# \EE_{\textup{LD}}(\PP) = O(n^2)$. \oprocend
\end{rmrk}

\subsection{Spatially-distributed functions, vector fields, and
  set-valued maps}

We are now in a position to discuss distributed control laws and
algorithms in formal terms.  From now on, we shall deal not only with
sets of distinct points, but also with tuples of elements of
$\real^2$.  Let~$\GG$ be a proximity graph function and let $Y$ be a
set.  A map $\map{f}{(\real^2)^n}{Y^n}$ is \emph{spatially distributed
  over $\GG$} if there exist maps $\map{\tilde f_i}{\real^2\times
  2^{(\real^2)^n}}{Y}$, $i\in\until{n}$, with the property that for
all $(p_1,\dots,p_n)\in(\real^2)^n$
\begin{equation*}
  f_i(p_1,\dots,p_n) = \tilde f_i(p_i, \setdef{p_j\in\real^2}{p_j
  \in\NN_{\GG(\{p_1,\dots,p_n\})}(p_i)} ) \, ,
\end{equation*}
where $f_i$ denotes the $i$th-component of~$f$.  A vector field $X$ on
$(\real^2)^n$ is \emph{spatially distributed over $\GG$} if its
associated map $\map{X}{(\real^2)^n}{(\real^2)^n}$, where the
canonical identification between the tangent space of $(\real^2)^n$
and $(\real^2)^n$ itself is understood, is spatially distributed in
the above sense.  Finally, a set-valued map
$\map{T}{(\real^2)^n}{2^{(\real^2)^n}}$ is \emph{spatially distributed
  over $\GG$} if there exist maps $\map{\tilde T_i}{\real^2\times
  2^{(\real^2)^n}}{2^{\real^2}}$, $i\in\until{n}$, with the property
that for all $(p_1,\dots,p_n)\in(\real^2)^n$
\begin{equation*}
  T_i(p_1,\dots,p_n) = \tilde T_i(p_i, \setdef{p_j\in\real^2}{p_j
  \in\NN_{\GG(\{p_1,\dots,p_n\})}(p_i)} ) \, ,
\end{equation*}
where $T_i$ denotes the $i$th-component of~$T$.

\begin{rmrk}
  In other words, to compute the $i$th component of a spatially-distributed
  function, vector field, or set-valued map at $(p_1,\dots,p_n)$, it is
  only required the knowledge of the vertex $p_i$ and the neighboring
  vertexes in the graph $\GG(\{p_1,\dots,p_n\})$. \oprocend
\end{rmrk}

We are now in a position to establish an important property of the
$r$-limited Delaunay graph.
\begin{lmm}
  Let $r\in\real_+$.  The map
  $\map{\NN_{\subscr{\GG}{LD}(\cdot,r)}}{(\real^2)^n}{\big[2^{(\real^2)^n}\big]^n}$,
  defined by
  \begin{gather*}
    (p_1,\dots,p_n) \mapsto
    (\NN_{\subscr{\GG}{LD}(\{p_1,\dots,p_n\},r)}(p_1), \dots,
    \NN_{\subscr{\GG}{LD}(\{p_1,\dots,p_n\},r)}(p_n)) \, ,
  \end{gather*}
  is spatially distributed over $\subscr{\GG}{disk}
  (\{p_1,\dots,p_n\},r)$.
\end{lmm}
\begin{proof}
  The result follows from
  Remark~\ref{re:jorge-is-smart}\ref{item-remark}.
\end{proof}
Loosely speaking, this lemma states that the $r$-limited Delaunay graph
$\subscr{\GG}{LD}$ can be computed in a spatially localized way: each agent
needs to know only the location of all other agents in a disk of radius
$r$.  This property is to be contrasted with the centralized computation
required to determine the $r$-Delaunay graph
$\subscr{\GG}{disk\intersect{D}}$. This requirement can be observed in
Figure~\ref{fig:counter-example}: if $p_j$ and $p_l$ are placed arbitrarily
close to the line joining $p_i$ and $p_k$, then, in order to decide if $p_l
\in \NN_{\subscr{\GG}{disk\intersect{D}}(\PP,r)}(p_i)$, in general it is
necessary to know the locations of all the other points in
$\{p_1,\dots,p_n\}$.

\subsection{Piecewise smooth sets and a generalized
  Conservation-of-Mass Law}

A set $S \subset \real^2$ is called \emph{strictly star-shaped} if
there exists a point $p\in S$ such that for all $s \in \partial S$ and
all $\lambda \in (0,1]$, one has that $\lambda p + (1-\lambda) s \in
\interior(S)$. A \emph{curve} $C$ in $\real^2$ is the image of a map
$\map{\gamma}{[a,b]}{\real^2}$. The map $\gamma$ is called a
\emph{parameterization of $C$}.  A curve
$\map{\gamma}{[a,b]}{\real^2}$ is \emph{simple} if it is not
self-intersecting, i.e., if $\gamma$ is injective on $(a,b)$.  A curve
is \emph{closed} if $\gamma(a) = \gamma (b)$. A set $\Omega \subset
\real^2$ is said to be \emph{piecewise smooth} if its boundary,
$\partial \Omega$, is a simple closed curve that admits a continuous
and piecewise smooth parameterization
$\map{\gamma}{\sphere^1}{\real^2}$.  Likewise, a collection of sets
$\setdef{\Omega(x) \subset \real^2}{x \in (a,b)}$ is said to be a
\emph{piecewise smooth family} if $\Omega(x)$ is piecewise smooth for
all $x \in (a,b)$, and there exists
$\map{\gamma}{\sphere^1\times(a,b)}{\real^2}$, $(\theta,x) \mapsto
\gamma(\theta,x)$, differentiable with respect to its second argument
such that for each $x \in (a,b)$, the map $\theta \mapsto
\gamma_x(\theta) =\gamma(\theta,x)$ is a continuous and piecewise
smooth parameterization of $\partial \Omega(x)$.  We refer to~$\gamma$
as a \emph{parameterization for the family $\setdef{\Omega(x) \subset
    \real^2}{x \in (a,b)}$}.

The following result is an extension of the integral form of the
Conservation-of-Mass Law in fluid mechanics~\cite{AJC-JEM:94}.
The proof is given in Appendix~\ref{sec:mass-conservation}.  Given a
curve $C$ parameterized by a piecewise smooth map $\gamma:[a,b]
\rightarrow C$, recall that the line integral of a function $f: C
\subset \real^2 \rightarrow \real$ over~$C$ is defined by
\begin{align*}
  \int_C f  
  = \int_{a}^b f(\gamma(t)) \, |\dot{\gamma}(t)| \, dt \, ,
\end{align*}
and it is independent of the selected parameterization.

\begin{prpstn}\label{prop:mass-conservation}
  Let $\{ \Omega(x) \subset{Q} \, | \, x\in (a,b)\}$ be a piecewise smooth
  family such that $\Omega(x)$ is strictly star-shaped for all $x \in
  (a,b)$.  Let the function $\phi : Q \times (a,b) \rightarrow \real$ be
  continuously differentiable with respect to its second argument for all
  $x \in (a,b)$ and almost all $q \in \Omega (x)$, and such that for each
  $x \in (a,b)$, the maps $q \mapsto \phi(q,x)$ and $q \mapsto
  \pder{\phi}{x}(q,x)$ are integrable on $\Omega(x)$.  Then, the function
  \begin{equation}\label{eq:function}
    (a,b) \ni x  \mapsto \int_{\Omega(x)} \phi(q,x) dq
  \end{equation}
  is continuously differentiable and
  \begin{equation*}
    \der{}{x} \int_{\Omega(x)} \phi(q,x) dq = \int_{\Omega(x)}
    \pder{\phi}{x}(q,x) dq 
    + \int_{\partial\Omega(x)} n^t \pder{\gamma}{x} \,
    \phi(\cdot,x) \, ,
  \end{equation*}
  where $n:\partial \Omega (x) \rightarrow \real^2$, $q \mapsto n(q)$,
  denotes the unit outward normal to $\partial\Omega(x)$ at $q \in
  \partial\Omega(x)$, and $\gamma: \sphere^1 \times (a,b) \rightarrow
  \real^2$ is a parameterization for the family $\setdef{\Omega(x)
    \subset \real^2}{x \in (a,b)}$.
\end{prpstn}

\begin{rmrk}\label{re:extending-proposition}
  Note that this result can be readily extended to any family of sets
  $\setdef{\Omega(x)}{x \in (a,b)}$ that admits a suitable
  decomposition into piecewise smooth families consisting of strictly
  star-shaped sets.  For instance, if $\setdef{\Omega_i(x)}{x \in
    (a,b)}$, $i \in \{ 1,2\}$ are two piecewise smooth families with
  strictly star-shaped sets and we consider the family $\Omega(x) =
  \Omega_1(x) \setminus \Omega_2(x)$, for all $x \in (a,b)$, then the
  same result holds for the function
  \begin{equation*}
    x \in (a,b) \mapsto 
    \int_{\Omega(x)}\phi(q,x) dq = \int_{\Omega_1(x)} \phi(q,x) dq -
    \int_{\Omega_2(x)} \phi(q,x) dq \,, 
  \end{equation*}
  by applying Proposition~\ref{prop:mass-conservation} to each summand
  on the right-hand side of the equality. \oprocend
\end{rmrk}

\section{limited-range locational optimization}
\label{sec:locational-optimization}

Let $Q$ be a simple convex polygon in $\real^2$ including its
interior.  The \emph{diameter of~$Q$} is defined as $\diam{Q} =
\max_{q,p \in Q} \| q-p\|$.  For $\delta,\eps\in\real_+$,
$\delta<\eps$, let $D_{[\delta,\eps]}(p) = \setdef{q \in
  \real^2}{\delta \le \|q-p\|\le \eps} $ denote the annulus in
$\real^2$ centered at $p$ of radius $\delta$ and $\eps$; it is also
convenient to define
$D_{[0,\eps]}(p)=B_\eps(p)=\setdef{q\in\real^2}{\|q-p\|\le\eps}$ and
$D_{[\delta,+\infty)}(p) = \setdef{q \in \real^2}{\delta \le \|q-p\|}
$.  Let $n_{B_\eps(p)}(q)$ denote the unit outward normal to
$B_\eps(p)$ at $q \in \partial B_\eps(p)$.  Given a set $S \subset Q$,
let $\indicator{S}$ denote the indicator function defined by
$\indicator{S} (q) = 1$ if $q \in S$, and $\indicator{S} (q)=0$ if $q
\not \in S$.

Throughout the rest of the paper, given a set of $n$ distinct points
$\PP=\{p_1,\dots,p_n\}$, we consider the restriction of the Voronoi
partition $\VV(\PP)$ generated by $\PP$ to the convex polygon $Q$,
$\{V_i(\PP) \cap Q \}_{i \in \until{n}}$. For ease of exposition, we
denote this restriction in the same way as the standard Voronoi
partition. Given a tuple $P=(p_1,\dots,p_n)\in Q^n$, recall that we
let $\PP = \{p_1,\dots,p_n\}$ denote the associated point set that
only contains the corresponding distinct points.

A \emph{density function} $\map{\phi}{Q}{\ov{\real}_+}$\ is a bounded
function on $Q$.  Given a set $S \subset Q$, let $\area{S}$ denote the
area of $S$ measured according to~$\phi$, i.e., $\area{S} = \int_{S}
\phi(q) dq$.  A \emph{performance function}
$\map{f}{\ov{\real}_+}{\real}$ is a non-increasing and piecewise
differentiable map with finite jump discontinuities at
$R_1,\dots,R_m\in\real_+$, with $R_1 < \dots <R_m$.  For convenience,
we set $R_0 = 0$ and $R_{m+1}=+\infty$, and write
\begin{align}
  \label{eq:function-f}
  f (x) = \sum_{\alpha=1}^{m+1} f_\alpha (x)
  \indicator{[R_{\alpha-1},R_\alpha)} (x) \, ,
\end{align}
where $\map{f_\alpha}{[R_{\alpha-1},R_\alpha]}{\real}$, $\alpha \in
\until{m+1}$ are non-increasing differentiable functions such that
$f_\alpha (R_\alpha) > f_{\alpha+1} (R_\alpha)$ for $\alpha \in
\until{m}$.  Given a density function $\phi$ and a performance
function $f$, we consider the \emph{multi-center function}
$\map{\HH}{Q^n}{\real}$ defined by
\begin{equation}
  \label{eq:HH}
  \HH (P) = \int_Q \max_{i\in \until{n}} f(\| q - p_i \| ) \phi(q) dq. 
\end{equation}
Note that $\HH$ is an aggregate objective function in the sense that it
depends on all the locations $p_1,\dots,p_n$.  It will be of interest
to find local maxima for $\HH$.  Note that the performance function
$f$ can be defined up to a constant $c \in \real$, since
\begin{align*}
  \int_Q \Big( \max_{i\in \until{n}} \big( f(\|q-p_i\|) + c\big) \Big)
  \phi(q) dq= \HH(P) + c \, \area{Q}\, ,
\end{align*}
and, therefore, this function and $\HH$ have the same local maxima.

\begin{rmrk}
  Maximizing the multi-center function is an optimal resource
  placement problem; it is interesting to draw an analogy with the
  optimal placement problem for large numbers of spatially-distributed
  sensors. In this setting, (1) $\HH$ provides the \emph{expected
    value of the sensing performance} provided by the group of sensors
  over any point in the environment~$Q$, where
  (2) the function $\phi$ is the \emph{distribution density function}
  representing a measure of information or probability that some event
  take place over~$Q$,
  and (3) $f$ describes the \emph{utility or sensing performance} of
  the sensors.  Because of noise and loss of resolution, the sensing
  performance at point~$q$ taken from $i$th sensor at the position
  $p_i$ degrades with the distance $\|q -p_i \|$ between $q$ and
  $p_i$. \oprocend
\end{rmrk}

Let us provide two equivalent expressions for the function $\HH$ over the
set $Q^n \setminus \Sc$, for $\Sc$ defined in equation~\eqref{eq:Sc}. Using
the definition of the Voronoi partition and the fact that $f$ is
non-increasing, $\HH$ can be rewritten as
\begin{align*}
  \HH (P) & = \sum_{i=1}^n \int_{V_i(P)} f(\| q - p_i \| ) \phi(q) dq
  \, , \quad P \in Q^n \setminus \Sc \, .
\end{align*}
Resorting to the expression of $f$ in~\eqref{eq:function-f}, we can
also rewrite $\HH$ as
\begin{align}
  \label{eq:HH-decomposed}
  \HH (P)
  &=\; \sum_{i=1}^n \sum_{\alpha=1}^{m+1} \int_{V_i(P) \cap
    D_{[R_{\alpha-1},R_\alpha]}(p_i)} f_\alpha (\| q - p_i \| )
  \phi(q) dq \, , \quad P \in Q^n \setminus \Sc \, .
\end{align}

We now analyze the smoothness properties of the multi-center function.
\begin{thrm} \label{the:smoothness-HH}
  Given a density function $\phi$ and a performance function $f$, the
  multi-center function $\HH$ is
  \begin{enumerate}
  \item globally Lipschitz on $Q^n$, and 
  \item continuously differentiable on $Q^n \setminus \Sc$, where for
    each $i \in \until{n}$
    \begin{align}
      \pder{\HH}{p_i}(P) = &\; \int_{V_i(P)} \pder{}{p_i} f (\| q - p_i
      \| ) \phi(q) dq  \nonumber \\
      & \; + \sum_{\alpha=1}^m \left( f_\alpha (R_\alpha ) - f_{\alpha+1}
        (R_\alpha ) \right) \Big( \sum_{k=1}^{M_i(2R_\alpha)}
      \int_{\arc_{i,k}(2R_\alpha)} n_{B_{R_\alpha}(p_i)}(q) \phi(q) dq
      \Big) \, ,
      \label{eq:partial-HH}
    \end{align}
    with $\arc_{i,k}(2R_\alpha)$, $k\in\until{M_i(2R_\alpha)}$ the
    arcs in the boundary of $V_i(P)\intersect B_{R_\alpha}(p_i)$.
  \end{enumerate}
\end{thrm}
\begin{proof}
  We start by proving fact (i).  Because $\max_{i\in \until{n}} \| q -
  p_i \| \leq \diam{Q}$ for all $q,p_1,\dots,p_n$ in $Q$, we can
  assume, without loss of generality, that $R_{m+1} = \diam{Q}$.
  Since the functions $f_\alpha$, $\alpha \in \until{m+1}$, are
  differentiable on $[R_{\alpha -1}, R_{\alpha}]$, they admit a
  non-increasing $C^1$-extension to $[0,R_{\alpha}]$, that we also
  denote by $f_\alpha$ for simplicity.  We then rewrite $\HH$ as
  \begin{equation*}
    \HH (P) = \sum_{\alpha=1}^{m+1} \int_Q f_{\alpha}(\dist(q,P)) \left(
      \indicator{[0,R_{\alpha})}(\dist(q,P))- \indicator{[0,R_{\alpha
          -1})}(\dist(q,P)\right)\phi(q) dq \,,
  \end{equation*}
  where $\dist(q,P) = \min_{i \in \until{n}}\| q -p_i\|$, for $P =
  (p_1, \dots p_n) \in Q^n$.  Since the finite sum of globally
  Lipschitz functions is globally Lipschitz, in what follows it
  suffices to prove that for $R \in [0,R_\alpha]$ and for $\alpha \in
  \until{m+1}$, the function
  \begin{equation*}
    \HH_{\alpha,R}(P) = \int_Q f_\alpha
    (\dist(q,P))\indicator{[0,R)}(\dist(q,P)) \phi(q)dq\,,
  \end{equation*}
  is globally Lipschitz.  To this end, we introduce a useful partition
  of $Q$.  For $S\subset Q$, recall $S^c=Q\setminus S$.  Given
  $P=(p_1,\dots,p_n)$, $P'=(p'_1,\dots,p'_n)$, define the following
  sets
  \begin{align*}
    S_1 & = \Big( \union\limits_{i \in \until{n}} B_R(p_i)\Big)
    \intersection \Big( \intersection\limits_{i \in
      \until{n}}B_R^c(p'_i) \Big)\,, \\
    S_2 &= \Big( \union\limits _{i \in \until{n}} B_R(p'_i)\Big)
    \intersection \Big( \intersection\limits_{i \in
      \until{n}}B_R^c(p_i) \Big) \,,\\
    S_3 &= \Big( \intersection\limits _{i \in \until{n}}B_R^c(p_i)
    \Big) \intersection \Big(
    \intersection\limits_{i \in \until{n}}B_R^c(p'_i)\Big)\,,\\
    S_4 &= \Big( \union\limits_{i \in \until{n}}B_R(p_i)\Big)
    \intersection \Big( \union\limits_{i \in \until{n}}
    B_R(p'_i)\Big)\,.
  \end{align*}
  Note that $S_1\union S_2 = (S_3\union S_4)^c$ and therefore $Q = S_1 \cup
  S_2 \cup S_3 \cup S_4$. Also, observe that $S_a \cap S_b = \emptyset$ for
  any $a$, $b\in \{1,2,3,4\}$, $a \neq b$. Accordingly, we write
  \begin{multline*}
    \HH_{\alpha,R}(P) -\HH_{\alpha,R}(P') \\
    = \sum_{a=1}^4 \int_{S_a} \left(
    f_{\alpha} (\dist(q,P)) \indicator{[0,R)}(\dist(q,P))
    -f_{\alpha }(\dist(q,P'))
    \indicator{[0,R)}(\dist(q,P'))\right)\phi(q)dq\,.
  \end{multline*}
  Now we upper bound each of the integrals in the above sum.  For $q
  \in S_3$, we have $\indicator{[0,R)}(\dist(q,P))=0$ and
  $\indicator{[0,R)}(\dist(q,P'))=0$, and therefore the integral over
  $S_3$ vanishes.  For $q \in S_4$, we have $\dist(q,P) \le R$ and
  $\dist(q,P') \le R$.  Thus,
  \begin{align*}
    \Big| \int_{S_4} \big( f_{\alpha} (\dist(q,P))
    &\indicator{[0,R)}(\dist(q,P)) - f_{\alpha }(\dist(q,P'))
    \indicator{[0,R)}(\dist(q,P'))\big)\phi(q)dq \; \Big|  \\
      & \le\; \int_{S_4}\, | \; f_{\alpha} (\dist(q,P)) -f_{\alpha
        }(\dist(q,P'))\, |\, \phi(q)dq \; \\
      & \le\; \left\|\frac{df_\alpha}{dx} \right\|_{[0,R_\alpha]}
      \int_{S_4}
      \,|\, \dist(q,P) - \dist(q,P') \,|\,\phi(q)dq \\
      &\le\; \left\|\frac{df_\alpha}{dx} \right\|_{[0,R_\alpha]} \|P -P' \|
      \int_{S_4} \phi(q) dq 
      \; \le \;\left\|\frac{df_\alpha}{dx} \right\|_{[0,R_\alpha]}
      \text{area}_{\phi}(Q)\;\|P -P' \|\,,
  \end{align*}
  where $ \| g \|_{[0,R_\alpha]}$ denotes the $L_{\infty}$-norm of
  $\map{g}{[0,R_\alpha]}{\real}$, and $ \| P - P'\|$ is the Euclidean
  norm of $P-P'$ as a vector in $\real^{2n}$.  Here we have made use
  of the fact that, for all $q\in Q$, the map $P\mapsto\dist(q,P)$ is
  globally Lipschitz with Lipschitz constant~$1$.  For $q \in S_1$, we
  have
  \begin{align*}
    \Big|\; \int_{S_1} \Big( f_{\alpha} (\dist(q,P))
    & \indicator{[0,R)}(\dist(q,P)) -  f_{\alpha }(\dist(q,P'))
    \indicator{[0,R)}(\dist(q,P'))\Big)\phi(q)dq \; \Big|  \\
    &\le\; \int_{S_1}\, |\,f_{\alpha} (\dist(q,P)) \,|\,\phi(q)dq
    \;\le\; \| \phi \|_{Q} \| f_\alpha \|_{[0,R_\alpha]} \int_{S_1} dq
    \\
    &\le\; \| \phi \|_{Q} \| f_\alpha \|_{[0,R_\alpha]} \sum_{i=1}^n
    \int_{B_R(p_i) \intersection ( \intersection_{j \in
        \until{n}}B_R^c(p'_j))} dq \\
    &\le\; \| \phi \|_{Q} \| f_\alpha \|_{[0,R_\alpha]} \sum_{i=1}^n
    \int_{B_R(p_i) \cap B_R^c(p'_i)} dq,
  \end{align*}
  where $\|\phi\|_{Q}=\max_{q\in Q} \phi(q)$.  Now, we observe that,
  for $\| p_i -p'_i \| \ge R$,
  \begin{equation}
    \int_{B_R(p_i)\cap B_R^c(p'_i)}dq \le \pi R^ 2  \le \pi \,
    \text{diam}(Q) \,\| p_i - p'_i \|\,. 
  \end{equation}
  On the other hand, for $\| p_i - p'_i\| \le R$,
  Lemma~\ref{le:overlap-area} in Appendix~\ref{app:overlap-area} shows
  that
  \begin{equation}
     \int_{B_R(p_i)\cap B_R^c(p'_i)}dq \le \tfrac{2\sqrt{3}+3}{3} R \,
     \| p_i -p'_i\|  \le  \tfrac{2\sqrt{3}+3}{3} \diam{Q} \, \| p_i
     -p'_i\|. 
  \end{equation}
  Therefore, since $\tfrac{2\sqrt{3}+3}{3}<\pi$, we have
  \begin{align*}
    \Big|\; \int_{S_1} \Big( f_{\alpha} & (\dist(q,P))
    \indicator{[0,R)}(\dist(q,P)) - f_{\alpha }(\dist(q,P'))
    \indicator{[0,R)}(\dist(q,P'))\Big)\phi(q)dq \; \Big|  \\
    &\le\; \pi \| \phi \|_{Q} \| f_\alpha \|_{[0,R_\alpha]} \diam{Q}
    \sum_{i=1}^n \| p_i-p'_i\| \;\le\; \frac{\pi}{\sqrt{n}} \| \phi \|_{Q}
    \| f_\alpha \|_{[0,R_\alpha]} \diam{Q} \; \| P -P'\| .
  \end{align*} 
  The integral over $S_2$ can be bounded in an analogous fashion.
  Summarizing, we have proved that $\HH_{\alpha,R}$ satisfies
  \begin{equation*}
    | \HH_{\alpha,R}(P) -\HH_{\alpha,R}(P') | \le L_{\alpha,R} \| P -
      P'\|\,, 
  \end{equation*}
  with $ L_{\alpha,R} = \frac{2\pi}{\sqrt{n}} \| \phi \|_Q \| f_\alpha
  \|_{[0, R_\alpha]} \diam{Q} + \left\|
    \frac{df_\alpha}{dx}\right\|_{[0,R_\alpha]}\text{area}_\phi (Q)$.  This
  concludes the proof of the statement that $\HH$ is globally Lipschitz on
  $Q^n$.
  
  Next, we prove fact (ii), that is, we prove that $\HH$ is
  continuously differentiable on $Q^n \setminus \Sc$ and we compute
  its partial derivative with respect to $p_i$.  Consider the
  expression~\eqref{eq:HH-decomposed} for the function~$\HH$.  Note
  that for each $i \in \until{n}$ and $\alpha \in \until{m+1}$, the
  function $(q,P) \mapsto f_\alpha (\| q - p_i \| )$ is continuously
  differentiable with respect to its second argument for all $P \in
  Q^n$ and almost all $q \in V_i(P) \cap
  D_{[R_{\alpha-1},R_\alpha]}(p_i)$. Note also that, for fixed $P \in
  Q^n$, both $q \mapsto f_\alpha (\| q - p_i \| )$ and $q \mapsto
  \pder{}{P} (f_\alpha (\| q - p_i \| ))$ are integrable on $V_i(P)
  \cap D_{[R_{\alpha-1},R_\alpha]}(p_i)$. Furthermore, if $P \not \in
  \Sc$, then the set
  \begin{multline*}
    \{q \in Q|\enspace \exists i,j \in \until{n} , \; i \neq j,  \;
      \text{such that} \\  
      \; \|q-p_i\| = \| q-p_j\| \le \|q-p_k\| \; \text{for} \; k \in
    \until{n} \setminus\{i,j\} \}
  \end{multline*}
  has measure zero. Therefore, $\setdef{V_i(P)}{P \in Q^n \setminus
    \Sc}$ is a piecewise smooth family for each $i \in \until{n}$.
  Since for each $\alpha \in \until{m+1}$, the balls
  $\setdef{B_{R_\alpha}(p_i)}{P \in Q^n }$ also define a piecewise
  smooth family, one concludes that the intersection $V_i \cap
  D_{[R_{\alpha-1},R_\alpha]}(p_i) = V_i \cap B_{R_\alpha}(p_i)
  \setminus V_i \cap B_{R_{\alpha-1}}(p_i)$, with $P \in Q^n \setminus
  \Sc$, can be written as the difference of two piecewise smooth
  families with strictly star-shaped sets. Applying now
  Proposition~\ref{prop:mass-conservation} (see also
  Remark~\ref{re:extending-proposition}), we deduce that each summand
  in equation~\eqref{eq:HH-decomposed} is continuously differentiable
  on $Q^n \setminus \Sc$. We now compute its partial derivative with
  respect to $p_i$, $i \in \until{n}$, as
  \begin{align*}
    \pder{\HH}{p_i}(P) = &\; \pder{}{p_i} \left( \sum_{\alpha=1}^{m+1}
      \int_{V_i(P) \cap D_{[R_{\alpha-1},R_\alpha]}(p_i)}
      f_\alpha (\| q - p_i \| ) \phi(q) dq \right) \\
    & \; + \pder{}{p_i} \left( \sum_{j \neq i} \sum_{\alpha=1}^{m+1}
      \int_{V_j(P) \cap D_{[R_{\alpha-1},R_\alpha]}(p_j)} f_\alpha (\| q - p_j
      \| ) \phi(q) dq \right) \, .
  \end{align*}
  For each $k \in \until{n}$ and each $\alpha \in \until{m+1}$, let
  $n_{k,\alpha}(q)$ denote the unit outward normal to $V_k(P) \cap
  B_{R_\alpha}(p_k)$ at $q$, and let $\gamma_{k,\alpha}:\sphere^1
  \times Q^n \setminus \Sc \rightarrow \real^2$ denote a
  parameterization for the family $\left\{V_k(P) \cap
    B_{R_\alpha}(p_k) \; | \; P \in Q^n \setminus \Sc \right\}$.
  Using Proposition~\ref{prop:mass-conservation}, the above expression
  is equal to
  \begin{align*}
    \pder{\HH}{p_i}(P) =&\; \sum_{\alpha=1}^{m+1} \int_{V_i(P) \cap
      D_{[R_{\alpha-1},R_\alpha]}(p_i)} \pder{}{p_i} f_\alpha (\| q -
    p_i \| ) \phi(q) dq \\
    &\;+ \sum_{\alpha=1}^{m+1} \int_{\partial \big( V_i(P) \cap
      B_{R_\alpha}(p_i) \big)} n_{i,\alpha}^t
    \pder{\gamma_{i,\alpha}}{p_i} f_\alpha (\dist (\cdot, p_i )) \,
    \phi  \\
    & \; - \sum_{\alpha=1}^{m+1} \int_{\partial \big( V_i(P) \cap
      B_{R_{\alpha-1}}(p_i) \big)} n_{i,\alpha-1}^t
    \pder{\gamma_{i,\alpha-1}}{p_i} f_\alpha (\dist (\cdot, p_i ))
    \, \phi \\
    &\; + \sum_{\alpha=1}^{m+1} \sum_{j \neq i} \int_{
      \begin{subarray}{c}
        \partial \big( V_j(P) \cap B_{R_\alpha}(p_j) \big) \cap
        \partial \big(V_i(P) \cap \, B_{R_\alpha}(p_i)\big)
      \end{subarray}   
      } n_{j,\alpha}^t \pder{\gamma_{j,\alpha}}{p_i} f_\alpha (\dist
    (\cdot, p_j )) \, \phi \\
    & \; - \sum_{\alpha=1}^{m+1} \sum_{j \neq i} \int_{
      \begin{subarray}{c}
        \partial \big( V_j(P) \cap B_{R_{\alpha-1}}(p_j) \big) \cap
        \partial \big(V_i(P) \cap \, B_{R_{\alpha-1}}(p_i)\big)
      \end{subarray}   
      } n_{j,\alpha-1}^t \pder{\gamma_{j,\alpha-1}}{p_i}
    f_\alpha(\dist (\cdot, p_j )) \, \phi \, ,
  \end{align*} 
  where recall that $\dist(q,p) = \|q-p\|$. For $\alpha \in
  \until{m+1}$, note that $\Delta_{ij}(2R_\alpha) = (V_i(P) \cap
  B_{R_\alpha}(p_i)) \cap (V_j(P) \cap B_{R_\alpha}(p_j)) \neq
  \emptyset$ if and only if $p_i$ and $p_j$ are neighbors according to
  the $2R_\alpha$-limited Delaunay graph
  $\subscr{\GG}{LD}(P,2R_\alpha)$.  In this case, there exist
  intervals $[\theta^-_{i,j}(P),\theta^+_{i,j}(P)]$ and
  $[\theta^-_{j,i}(P),\theta^+_{j,i}(P)]$ depending smoothly on $P$
  over an open set of $Q^n \setminus \Sc$ such that
  \begin{align*}
    \theta \in [\theta^-_{i,j}(P),\theta^+_{i,j}(P)] \mapsto
    \gamma_{i,\alpha}(\theta,P) \, , \quad \ov{\theta} \in
    [\theta^-_{j,i}(P),\theta^+_{j,i}(P)] \mapsto
    \gamma_{j,\alpha}(\ov{\theta},P) \, ,
  \end{align*}
  are two parameterizations of the set $(V_i(P) \cap
  B_{R_\alpha}(p_i)) \cap (V_j(P) \cap B_{R_\alpha}(p_j))$. Resorting
  to the implicit function theorem, one can show that there exists a
  function $h: \sphere^1 \times U \rightarrow \sphere^1$,
  $h([\theta^-_{j,i}(P),\theta^+_{j,i}(P)],P) =
  [\theta^-_{i,j}(P),\theta^+_{i,j}(P)]$, such that
  $\gamma_{j,\alpha}(\theta,P) = \gamma_{i,\alpha}(h(\theta,P),P)$ for
  $\theta \in [\theta^-_{j,i}(P),\theta^+_{j,i}(P)]$. From here, we
  deduce that $n_{j,\alpha}^t \pder{\gamma_{j,\alpha}}{p_i} =
  n_{j,\alpha}^t \left( \pder{\gamma_{i,\alpha}}{\theta} \pder{h}{p_i}
    + \pder{\gamma_{i,\alpha}}{p_i} \right) = n_{j,\alpha}^t
  \pder{\gamma_{i,\alpha}}{p_i}$, since $n_{j,\alpha}$ and $
  \pder{\gamma_{i,\alpha}}{\theta}$ are orthogonal.  Therefore, if
  $p_j \in \NN_{\subscr{\GG}{LD}(P,2R_\alpha)}(p_i)$, we have
  \begin{multline*}
    \int_{
      \begin{subarray}{c}
        \partial \big( V_j(P) \cap B_{R_\alpha}(p_j) \cap V_i(P) \cap
        B_{R_\alpha}(p_i) \big)
      \end{subarray}
      } n_{j,\alpha}^t \pder{\gamma_{j,\alpha}}{p_i} f_\alpha(\dist
    (\cdot, p_j )) \, \phi \\
    = - \int_{
      \begin{subarray}{c}
        \partial \big( V_i(P) \cap B_{R_\alpha}(p_i) \cap V_j(P) \cap
        B_{R_\alpha}(p_j) \big)
      \end{subarray}
      } n_{i,\alpha}^t \pder{\gamma_{i,\alpha}}{p_i}
    f_\alpha(\dist(\cdot,p_i)) \, \phi \, ,
  \end{multline*}
  since $n_{i,\alpha}(q)=-n_{j,\alpha}(q)$ and $\| q - p_i \| = \| q -
  p_j \|$ for all $q \in \partial \big( V_i(P) \cap B_{R_\alpha}(p_i)
  \cap V_j(P) \cap B_{R_\alpha}(p_j) \big)$.  Moreover, notice that if
  $p_i$ moves, the motion \textemdash{}projected to the normal
  direction $n_{i,\alpha}$\textemdash{} of the points in the arcs
  $\{\arc_{i,1}(2R_\alpha), \dots, \arc_{i,M_i(2 R_\alpha)}(2R_\alpha)
  \} \subset \partial (V_i(P) \cap B_{R_\alpha}(p_i))$ is exactly the
  same as the motion of $p_i$, i.e., $n_{i,\alpha}^t
  \pder{\gamma_{i,\alpha}}{p_i} = n_{i,\alpha}^t$ over
  $\arc_{i,1}(2R_\alpha) \cup \dots \cup \arc_{i,M_i(2
    R_\alpha)}(2R_\alpha)$.  Using this fact, the expression for the
  partial derivative of $\HH$ with respect to $p_i$ can be rewritten
  as
  \begin{multline*}
    \pder{\HH}{p_i}(P) = \int_{V_i(P)} \pder{}{p_i} f (\| q - p_i \| )
    \phi(q) dq + \sum_{\alpha=1}^{m+1} \Big( \sum_{k=1}^{M_i(2
      R_\alpha)} \int_{\arc_{i,k}(2R_\alpha)}
    n_{B_{R_{\alpha}}(p_i)} f_\alpha (R_\alpha) \, \phi  \\
    - \sum_{k=1}^{M_i(2 R_{\alpha-1})}
    \int_{\arc_{i,k}(2R_{\alpha-1})} n_{B_{R_{\alpha-1}}(p_i)}
    f_\alpha (R_{\alpha-1}) \, \phi \Big) .
  \end{multline*}     
  The final result is a rearrangement of the terms in this equation.
\end{proof}

\begin{rmrk}
  For a constant density function, $q\mapsto \phi (q) = c \in
  \real_+$, each line integral
  \begin{align*}
    \int_{\arc(2R)} n_{B_{R}(p)} \, \phi
  \end{align*}
  computed over the $\arc(2R)$ described by $[\theta_-,\theta_+]\ni
  \theta \mapsto p+R (\cos \theta,\sin \theta) \in \real^2$, equals
  \begin{equation*}
    c R \int_{\theta_-}^{\theta_+} (\cos \theta, \sin \theta) d\theta = 2 c
    R \sin\!\Big(\frac{\theta_+-\theta_-}{2}\Big) 
    \Big(
    \cos\!\Big(\frac{\theta_+ +\theta_-}{2}\Big),
    \sin\!\Big(\frac{\theta_+ +\theta_-}{2}\Big) \Big) . \eqoprocend
  \end{equation*}
\end{rmrk}

For particular choices of performance function, the corresponding
multi-center function and its gradient have different features.  We here
explore some interesting cases:
\begin{description}
\item[Centroid problem] If the performance function $f$ is piecewise
  differentiable with no jump discontinuities, then all the terms in
  the second summand of equation~\eqref{eq:partial-HH} vanish and one
  obtains
  \begin{align*}
    \pder{\HH}{p_i}(P) = \int_{V_i(P)} \pder{}{p_i} f (\| q - p_i \| )
    \phi(q) dq \, .
  \end{align*}
  This is the result known in the locational optimization
  literature~\cite{YA:91,AO-BB-KS-SNC:00,QD-VF-MG:99}. In particular,
  if $f(x) = -x^2$, the multi-center function $\HH$ reads
  \begin{align*}
    \HH (P) & = - \sum_{i=1}^n \int_{V_i(P)} \| q - p_i \|^2 \phi(q)
    dq \eqdef - \sum_{i=1}^n J_{V_i,p_i},
  \end{align*}
  where $J_{W,p}$ denotes the polar moment of inertia of the set $W
  \subset Q$ about the point $p$. Additionally, the gradient of~$\HH$
  is
  \begin{align*}
    \pder{\HH}{p_i}(P) & = 2 \int_{V_i(P)} (q-p_i) \phi(q) dq = 2
    \MM_{V_i(P)} (\CM_{V_i(P)} - p_i) \, .
  \end{align*}
  Here $\MM_W$ and $\CM_W$ denote, respectively, the mass and the
  center of mass with respect to the density function $\phi$ of the
  set $W \subset Q$. The critical points of $\HH$ are configurations
  $P \in Q^n$ such that $p_i = \CM_{V_i(P)}$ for all $i\in \until{n}$.
  Such configurations are called \emph{centroidal Voronoi
    configurations}, see~\cite{QD-VF-MG:99}.
      
\item[Area problem] On the other hand, if one takes $f(x) =
  \indicator{[0,R]}(x)$, the indicator function of the set $[0,R]$,
  then the multi-center function $\HH$ corresponds to the area,
  measured according to~$\phi$, covered by the union of the $n$ balls
  $B_R(p_1),\dots,B_R(p_n)$, that is,
  \begin{align*}
    \HH (P) = \area{\cup_{i=1}^n B_{R}(p_i)} \, .
  \end{align*}
  In this case, the first term in equation~\eqref{eq:partial-HH}
  vanishes and one obtains
  \begin{align*}
    \pder{\HH}{p_i}(P) = \sum_{k=1}^{M_i(2 R)} \int_{\arc_{i,k}(2R)}
    n_{B_{R}(p_i)} \, \phi \, .
  \end{align*}
  Given a configuration $P\in Q^n$, if the $i$th agent is surrounded
  by neighbors in the graph $\subscr{\GG}{LD}(P,2R)$ in such a way
  that $M_i(2R)=0$, then the multi-center function $\HH$ does not
  depend on $p_i$.  This situation is depicted in
  Figure~\ref{fig:limited-Delaunay} (see example on the right) and
  captures the fact that the total area covered by the agents is not
  affected by an infinitesimal displacement of the $i$th agent.
      
\item[Mixed centroid-area problem] Consider the case when the function
  $f$ is given by $x \mapsto - x^2 \, \indicator{[0,R)}(x) + b \cdot
  \indicator{[R,+\infty)}(x)$, for $b \le -R^2$.  The multi-center
  function then takes the form
  \begin{align*}
    \HH (P) &= - \sum_{i=1}^n J_{V_i(P) \cap B_R(p_i),p_i} + b \,
    \area{Q\setminus\cup_{i=1}^n B_{R}(p_i)} \, ,
  \end{align*}
  and its partial derivative with respect to the position of the $i$th
  agent is
  \begin{equation*}
    \pder{\HH}{p_i} (P) = 2 \MM_{V_i(P) \cap B_R(p_i)} (\CM_{V_i(P) \cap
      B_R(p_i)} - p_i) - (R^2+b) \! \sum_{k=1}^{M_i(2 R)}
    \hspace*{-3pt} \int_{\arc_{i,k}(2R)} \!  n_{B_{R}(p_i)} \, 
    \phi \, .
  \end{equation*}
  In the particular case when $b=-R^2$, the function $x \mapsto
  f(x) = -x^2 \indicator{[0,R)}(x) - R^2 \cdot
  \indicator{[R,+\infty)}(x)$ is continuous and therefore the gradient
  of~$\HH$ takes the form
  \begin{align*}
    \pder{\HH}{p_i}(P) & = 2 \int_{V_i(P) \cap B_R(p_i)} (q-p_i)
    \phi(q) dq = 2 \MM_{V_i(P) \cap B_R(p_i)} (\CM_{V_i(P) \cap
      B_R(p_i)} - p_i) \, .
  \end{align*}
  Note that, in this case, the critical points of $\HH$ are configurations
  $P \in Q^n$ such that $p_i = \CM_{V_i(P) \cap B_R(p_i)}$ for all $i\in
  \until{n}$. We refer to such configurations as \emph{$R$-centroidal
    Voronoi configurations}.  For $R \ge \diam{Q}$, $R$-centroidal Voronoi
  configurations coincide with the standard centroidal Voronoi
  configurations over~$Q$.
\end{description}

We can now characterize the results in Theorem~\ref{the:smoothness-HH}
in terms of the notion of spatially-distributed computations
introduced in Section~\ref{sec:geographs}.

\begin{crllr}\label{co:Delaunay-distributed}
  Let $\phi$ and $f$ be a density and a performance function,
  respectively.  The gradient of~$\HH$ with respect to the agents'
  location $P\in Q^n$ is spatially distributed over the Delaunay graph
  $\subscr{\GG}{D}(P)$.  Furthermore, if $f(x) = b$ for all $x \ge R$,
  then the gradient of $\HH$ with respect to the agents' location is
  spatially distributed over the $2R$-limited Delaunay graph
  $\subscr{\GG}{LD}(P,2R)$.
\end{crllr}
\begin{proof}
  In general, the partial derivative of $\HH$ with respect to the
  $i$th agent (cf. equation~\eqref{eq:partial-HH}) depends on the
  position $p_i$ and on the position of all neighbors of $p_i$ in the
  graph $\subscr{\GG}{D}$. If, in addition, $f(x) = b$, for all $x \ge
  R$, then necessarily $R_\alpha < R$, $\alpha \in \until{m}$, and
  \begin{equation*}
    \int_{V_i(P)} \pder{}{p_i} f (\| q - p_i \| ) \phi(q) dq 
    = \int_{V_i(P) \cap B_R(p_i)} \pder{}{p_i} f (\| q - p_i \| )
    \phi(q)  dq \,.
  \end{equation*}
  Therefore, the expression for ${\partial \HH}/{\partial p_i}$ in
  equation~\eqref{eq:partial-HH} can be computed with the knowledge
  $p_i$ and of its neighbors in the graph $\subscr{\GG}{LD}(P,2R)$.
\end{proof}

This corollary states that information about all neighbors in
$\subscr{\GG}{D}$ is required for objective functions $\HH$ corresponding
to arbitrary performance functions $f$. In the next proposition we
explore what can be done with only information about the neighbors in
the $2R$-limited Delaunay graph $\subscr{\GG}{LD}(2R)$.

\begin{prpstn}\label{prop:approximation-HH} 
  Let $f$ be a performance function and, without loss of generality,
  assume $f(0)=0$.  For $r \in ]0, 2 \diam{Q}]$, define the
  performance function $\map{f_{\frac{r}{2}}}{\ov{\real}_+}{\real}$
  given by $f_{\frac{r}{2}}(x) = f(x)$ for $x < \frac{r}{2}$ and
  $f_{\frac{r}{2}}(x) = f (\diam{Q})$ for $x \ge \frac{r}{2}$.  Let
  $\HH_{\frac{r}{2}}$ be the multi-center function associated to the
  performance function $f_{\frac{r}{2}}$.  Then, for all $P \in Q^n$,
  \begin{subequations}\label{eq:sandwichs}
    \begin{align}
      & \HH_{\frac{r}{2}} (P) \le \HH (P) \le \beta \,
      \HH_{\frac{r}{2}} (P) < 0 \, ,
      \label{eq:sandwich-multiplicative}
      \\
      & \HH_{\frac{r}{2}} (P) \le \HH (P) \le \HH_{\frac{r}{2}} (P) +
      \Pi(P) \, ,
      \label{eq:sandwich-additive}
    \end{align}
  \end{subequations}
  where $\beta = \frac{f(\frac{r}{2})}{f (\diam{Q})} \in [0,1]$ and
  $\Pi:Q^n \rightarrow [0,\kappa] \subset \real$, $\Pi(P) =
  (f(\frac{r}{2})- f (\diam{Q})) \area{Q \setminus \cup_{i=1}^n
    B_{\frac{r}{2}}(p_i)}$, with $\kappa = (f(\frac{r}{2})- f
  (\diam{Q})) \area{Q}$.
\end{prpstn}
\begin{proof}
  Clearly, $f_{\frac{r}{2}}$ is a performance function as it is
  non-increasing and piecewise differentiable with finite jump
  discontinuities.  Let $b = f (\diam{Q})$ and note that $f(x) \ge b$
  for all $x \in [0,\diam{Q}]$.  By construction, it is clear that
  $f_{\frac{r}{2}}(x) \le f(x)$ for all $x \in [0,\diam{Q}]$.  Since
  $\| q -p \| \le \diam{Q}$ for all $q, p \in Q$, we conclude that
  $\HH_{\frac{r}{2}} (P) \le \HH (P)$.  Now, consider the function
  $\tilde{f}(x) = \beta f_{\frac{r}{2}}(x)$. Note that $\tilde{f}(x) =
  \beta f(x) \ge f(x)$ for $x < \frac{r}{2}$, and $\tilde{f}(x) =
  \beta b = f(\frac{r}{2}) \ge f(x)$ for $x \ge \frac{r}{2}$.
  Therefore,
  \begin{align*}
    \beta \, \HH_{\frac{r}{2}}(P) & = \int_Q \max_{i\in \until{n}}
    \tilde{f}(\|q-p_i\|) \phi(q)dq \ge \int_Q \max_{i\in \until{n}}
    f(\|q-p_i\|) \phi(q)dq = \HH (P) \, ,
  \end{align*}
  which concludes the proof of the first chain of inequalities.  To
  prove the second chain of inequalities, consider the difference
  \begin{equation*}
    \HH (P) - \HH_{\frac{r}{2}} (P) 
    = \sum_{i=1}^n \int_{V_i(P) \cap (Q \setminus
    B_{\frac{r}{2}}(p_i))} (f(\|q-p_i\|) - b) \phi(q) dq \, .
  \end{equation*}
  For $q \in V_i(P) \cap (Q \setminus B_{\frac{r}{2}}(p_i))$, the
  non-increasing property of $f$ implies that $f(\|q-p_i\|) - b \le
  f(\frac{r}{2}) - b$.  Therefore,
  \begin{equation*}
    \HH (P) - \HH_{\frac{r}{2}} (P) \le \sum_{i=1}^n \int_{V_i(P) \cap (Q
      \setminus B_{\frac{r}{2}}(p_i))} ( f(\tfrac{r}{2}) - b) \phi(q) dq 
    = \int_{Q \setminus \cup_{i=1}^n B_{\frac{r}{2}}(p_i)} (
      f(\tfrac{r}{2}) - b) \phi(q) dq = \Pi (P) \, .
  \end{equation*}
\end{proof}

\begin{rmrk}
  The inequalities in~\eqref{eq:sandwichs} provide, respectively,
  constant-factor and additive approximations of the value of the
  multi-center function $\HH$ by the value of the function
  $\HH_{\frac{r}{2}}$. These approximations will play an important role in
  Section~\ref{sec:spat-dist-algorithms} when we discuss the continuous and
  discrete-time implementations of spatially-distributed coordination
  algorithms.  \oprocend
\end{rmrk}

The next result provides one more useful indication of the
relationship between multi-center functions associated to certain
performance functions.
\begin{prpstn}\label{prop:nice-observation}
  Let $f$ and $f_{\frac{r}{2}}$ be performance functions, and $\HH$
  and $\HH_{\frac{r}{2}}$ be the corresponding multi-center functions,
  defined as in Proposition~\ref{prop:approximation-HH}.  Let $P^* =
  (p_1^*,\dots,p_n^*) \in Q^n$ be a local maximum of
  $\HH_{\frac{r}{2}}$ such that $Q \subset \cup_{i \in \until{n}}
  B_{\frac{r}{2}}(p^*_i)$. Then $\HH(P^*) = \HH_{\frac{r}{2}} (P^*)$
  and $P^*$ is a local maximum of the aggregate objective function~$\HH$.
\end{prpstn}
\begin{proof}
  If $Q \subset \cup_{i \in \until{n}} B_{\frac{r}{2}}(p^*_i)$, then
  from equation~\eqref{eq:sandwich-additive} we deduce that $\HH(P^*)
  = \HH_{\frac{r}{2}} (P^*)$. Moreover, one can also show that $V_i(P)
  \subset B_{\frac{r}{2}}(p^*_i)$ for all $i \in \until{n}$, and
  therefore $V_i(P) = V_i(P) \cap B_{\frac{r}{2}}(p^*_i)$. As a
  consequence, the $r$-limited Delaunay graph
  $\subscr{\GG}{LD}(P^*,r)$ and the Delaunay graph
  $\subscr{\GG}{D}(P^*)$ coincide, and the gradients of both $\HH$ and
  $\HH_{\frac{r}{2}}$ vanish at~$P^*$.
\end{proof}

The importance of Proposition~\ref{prop:nice-observation} lies in the
fact that, by following the gradient of the
function~$\HH_{\frac{r}{2}}$ (where, along the evolution, the
inclusion $Q \subset \cup_{i \in \until{n}} B_{\frac{r}{2}}(p_i)$ may
not be verified and each agent only operates with the knowledge of (i)
the positions of other agents up to a distance $r$ of its own
position, and (ii) the events taking place at up to distance
$\frac{r}{2}$ of its own position), the agents may eventually find a
local maximum of the original multi-center function~$\HH$.

We end this section by presenting a useful result in the 1-center
case, i.e., when there is a single agent ($n=1$).  For a convex
polygon $W$, define the function $\HH_1(p,W) = \int_{W} f(\| q - p \|)
\phi(q) dq$.  The following lemma proves that the points in the
boundary of $W$ are not local maxima of $\HH_1 (\cdot,W): W
\rightarrow \real$.

\begin{lmm}\label{le:maxima-not-in-boundary}
  Let $W$ be a convex polygon, and consider the function $\HH_1
  (\cdot,W): W \rightarrow \real$. Let $p_0 \in \partial W$. Then the
  gradient of $\HH_1$ at $p_0$ is non-vanishing $\pder{\HH_1
    (\cdot,W)}{p} (p_0) \neq 0$, and points toward $\interior (W)$.
\end{lmm}
\begin{proof}
  The function $p \mapsto \HH_1(p,W)$ is differentiable over $W$, and
  its derivative is given by
  \begin{align*}
    \pder{\HH_1 (\cdot,W)}{p} & = \int_{W} \pder{}{p} f (\| q - p \| )
    \phi(q) dq + \sum_{\alpha=1}^{m+1} \int_{\partial (W \cap
      B_{R_\alpha}(p))} n_\alpha^t \pder{\gamma_{\alpha}}{p}
    f_\alpha(\dist(\cdot,p)) \; \phi \\
    & \; - \sum_{\alpha=1}^{m+1} \int_{\partial (W \cap
      B_{R_{\alpha-1}(p)})} n_{\alpha-1}^t \pder{\gamma_{\alpha-1}}{p}
    f_\alpha(\dist(\cdot,p)) \; \phi \, .
  \end{align*}
  Let $M(R_\alpha) \ge 0$ denote the number of distinct arcs
  $\arc_{1}(2R_\alpha),\dots,\arc_{M(R_\alpha)}(2R_\alpha)$ of radius
  $R_\alpha$ in $\partial (W \cap B_{R_\alpha}(p))$. After some
  simplifications, we rewrite the expression for the gradient at $p_0$
  as
  \begin{align}\label{eq:lemma-auxiliary}
    - \int_{W} f'(\| q - p_0 \| ) \frac{q-p_0}{\|q-p_0\|}\phi(q) dq +
    \sum_{\alpha=1}^{m} (f_\alpha (R_\alpha) - f_{\alpha+1}
    (R_\alpha)) \sum_{l=1}^{M(R_\alpha)} \int_{\arc_{l}(2R_\alpha)}
    n_{\alpha} \, \phi \, ,
  \end{align}
  where $n_\alpha$ denotes the outward normal to $B_{R_\alpha}(p_0)$.
  Since $W$ is convex, it is defined as the intersection of some
  hyperplanes $H_1, \dots, H_d$ via the equations $H_\zeta (q) =
  A_\zeta q + b_\zeta \ge 0$, where $A_\zeta$ is a
  $2\!\times\!2$-matrix and $b_\zeta \in \real$, for $\zeta \in
  \until{d}$.  To show that $\pder{\HH_1 (\cdot,W)}{p} (p_0) \neq 0$
  and points toward $\interior (W)$, we consider its inner product
  with the direction given by the each line $A_{\zeta_*} q +
  b_{\zeta_*}=0$ such that $H_{\zeta_*}(p_0)=0$. Let us therefore
  consider
  \begin{align*}
    A_{\zeta_*} \left( - \int_{W} f'(\| q - p_0 \| )
      \frac{q-p_0}{\|q-p_0\|}\phi(q) dq + \sum_{\alpha=1}^{m}
      (f_\alpha (R_\alpha) - f_{\alpha+1} (R_\alpha))
      \sum_{l=1}^{M(R_\alpha)} \int_{\arc_{l}(2R_\alpha)} n_{\alpha}
      \; \phi \right)  \\
    = - \int_{W} f'(\| q - p_0 \| ) \frac{A_{\zeta_*} q+
      b_{\zeta_*}}{\|q-p_0\|}\phi(q) dq + \sum_{\alpha=1}^{m}
    (f_\alpha (R_\alpha) - f_{\alpha+1} (R_\alpha))
    \sum_{l=1}^{M(R_\alpha)} \int_{\arc_{l}(2R_\alpha)}
    \frac{A_{\zeta_*} (\cdot) + b_{\zeta_*}}{\dist(\cdot,p_0)} \; \phi
    \, ,
  \end{align*}
  where we have used the fact that $n_\alpha (q) = (q-p_0)/\| q-p_0\|$
  for each $q \in \partial (W \cap B_{R_\alpha}(p_0))$.  Since the
  function~$f$ is non-increasing, then its derivative is negative
  almost everywhere, and the jump discontinuities $f_\alpha (R_\alpha)
  - f_{\alpha+1} (R_\alpha)$ are positive for all $\alpha \in
  \until{m}$.  Finally, note that $A_{\zeta_*} q+ b_{\zeta_*} >0$ in
  the interior of~$W$.  Therefore, we conclude that $A_{\zeta_*}
  \left( \pder{\HH_1 (\cdot,W)}{p} (p_0)\right) >0$ for all $\zeta_*$
  such that $H_{\zeta_*}(p_0) = 0$, i.e., $\pder{\HH_1 (\cdot,W)}{p}
  (p_0) \neq 0$ and points toward $\interior (W)$.
\end{proof}

\section{Design of spatially-distributed algorithms for coverage control}
\label{sec:spat-dist-algorithms}

In this section, we develop continuous and discrete-time implementations of
the gradient ascent for a general aggregate objective function~$\HH$.

\subsection{Continuous-time implementations}\label{se:continuous-time} 
Assume the agents' location obeys a first order dynamical behavior
described by
\begin{equation*}
  \dot{p}_i = u_i.
\end{equation*}
Consider $\HH$ an aggregate objective function to be maximized and impose
that the location $p_i$ follows the gradient ascent given
by~\eqref{eq:partial-HH}.  In more precise terms, we set up the
following control law defined over the set $Q^n \setminus \Sc$
\begin{align} \label{eq:continuous-lloyd}
  u_i = \pder{\HH}{p_i} (P) \, ,
\end{align}
where we assume that the partition $\VV(P)=\{V_1,\dots,V_n\}$ is
continuously updated. One can prove the following result.

\begin{prpstn}[Continuous-time Lloyd ascent]
  Consider the gradient vector field on $Q^n \setminus \Sc$ defined by
  equation~\eqref{eq:continuous-lloyd}. Then
  \begin{enumerate}
  \item For a general performance function $f$, the gradient vector
    field is spatially distributed over the Delaunay
    graph~$\subscr{\GG}{D}(\PP)$.  If, in addition, the performance
    function verifies $f(x) = b$ for all $x \ge R$, then the vector
    field is spatially distributed over the $2R$-limited Delaunay
    graph $\subscr{\GG}{LD}(P,2R)$.
  \item The agents' location evolving
    under~\eqref{eq:continuous-lloyd} starting at $P_0 \in Q^n
    \setminus \Sc$ remains in $Q^n \setminus \Sc$ and converges
    asymptotically to the set of critical points of the aggregate objective
    function~$\HH$.  Assuming this set is finite, the agents' location
    converges to a critical point of~$\HH$.
\end{enumerate}
\end{prpstn}

\begin{proof} 
  The statement (i) is a transcription of
  Corollary~\ref{co:Delaunay-distributed}. To prove the statement
  (ii), let $t \in \ov{\real}_+ \mapsto P(t) \in Q^n$ denote the
  solution to the initial value problem $\dot{p}_i =
  \pder{\HH_{\frac{r}{2}}}{p_i} (P)$, $i \in \until{n}$, $P(0)=P_0$.
  We reason by contradiction.  Assume that there exists $t_* \in
  \ov{\real}_+$ and $i,j\in \until{n}$ such that $p_i(t_*) =
  p_j(t_*)$, i.e., $P(t_*) \in \Sc$.  Let $v$ be the direction given
  by $v = \lim_{t \rightarrow t_*} \frac{p_i(t)-p_j(t)}{\|
    p_i(t)-p_j(t)\|}$.  Let $\eps>0$ sufficiently small such that, for
  all $t \in ]t_*-\eps,t_*[$, $p_i(t)$ and $p_j(t)$ are neighbors in
  the graph $\subscr{\GG}{LD}(P(t),r)$.  Then one can show that
  \begin{align}\label{eq:continuous-not-in-S}
    v \cdot \lim_{t \rightarrow t_*} \pder{\HH_{\frac{r}{2}}}{p_i}
    (P(t)) > 0 \, , \quad%
    v \cdot \lim_{t \rightarrow t_*} \pder{\HH_{\frac{r}{2}}}{p_j}
    (P(t)) < 0 \, .
  \end{align}
  Indeed, if $n$ denotes the orthogonal line to $v$, and $H_{i,n}$ and
  $H_{j,n}$ denote the associated hyperplanes having $v$ pointing
  inward and outward respectively, then, reasoning as in the proof of
  Lemma~\ref{le:maxima-not-in-boundary}, one proves that $\lim_{t
    \rightarrow t_*} \pder{\HH_{\frac{r}{2}}}{p_i} (P(t))$ points
  toward $\interior (V_i(P(t_*)) \cap B_{\frac{r}{2}}(p_i(t_*)) \cap
  H_{i,n} )$, and $\lim_{t \rightarrow t_*}
  \pder{\HH_{\frac{r}{2}}}{p_j} (P(t))$ points toward $\interior
  (V_j(P(t_*)) \cap B_{\frac{r}{2}}(p_j(t_*)) \cap H_{j,n} )$.  From
  equation~\eqref{eq:continuous-not-in-S}, we deduce that for all $t$
  sufficiently close to $t_*$, we have $(p_i(t)-p_j(t))\cdot \left(
    \dot{p}_i (t) - \dot{p}_j (t)\right) > 0$, which contradicts
  $P(t_*) \in \Sc$.  One can resort to a similar argument to guarantee
  that there is no configuration belonging to $\Sc$ in the
  $\omega$-limit set of the curve $t \mapsto P(t)$. The convergence
  result to the set of critical points of $\HH_{\frac{r}{2}}$ is an
  application of LaSalle Invariance Principle~\cite{HKK:96}.
\end{proof}

\begin{rmrk}
  Note that this gradient ascent is not guaranteed to find the global
  maximum.  For example, in the vector quantization and signal processing
  literature~\cite{RMG-DLN:98}, it is known that for ``bimodal''
  distribution density functions, the solution to the gradient flow reaches
  local maxima where the number of agents allocated to the two region of
  maxima are not optimally partitioned. \oprocend
\end{rmrk}

In a practical setting, the sensing and/or communication capabilities
of a network agent are restricted to a bounded region specified by a
finite radius $r>0$. Therefore, instead of maximizing the multi-center
function $\HH$, we set up the continuous-time algorithm given by
equation~\eqref{eq:continuous-lloyd} with the function
$\HH_{\frac{r}{2}}$. This latter algorithm has the advantage of being
spatially distributed over the $r$-limited Delaunay graph
$\subscr{\GG}{LD}(P,r)$, and providing an approximation of the
behavior for the multi-center function~$\HH$ (cf.
Proposition~\ref{prop:approximation-HH}).

\subsection{Discrete-time implementations}\label{se:discrete-time}

We start by reviewing some notions on discrete-time algorithms
following~\cite{DGL:84}. An \emph{algorithm on $Q^n$} is a set-valued
map $T : Q^n \rightarrow 2^{Q^n}$. Note that a map from $Q^n$ to $Q^n$
can be interpreted as a singleton-valued map.  For any initial $P_0
\in Q^n$, an algorithm $T$ generates feasible sequences of
configurations in the following way: given $P_n \in Q^n$, the map $T$
yields $T(P_n) \subset Q^n$. From this set, an arbitrary element
$P_{n+1}$ may be selected.  In other words,
\begin{align}\label{eq:algorithm}
  P_{n+1} \in T (P_n) \, , \quad n \in \natural \cup \{0 \} \,.
\end{align}
An algorithm $T$ is said to be \emph{closed at $P\in Q^n$} if for all
convergent sequences $P_k \rightarrow P$, $P'_k \rightarrow P'$ such
that $P'_k \in T(P_k)$, one has that $P' \in T(P)$. An algorithm is
said to be \emph{closed on $W \subset Q^n$} if it is closed at $P$,
for all $P \in W$. In particular, every continuous map $T:Q^n
\rightarrow Q^n$ is closed on~$Q^n$.
A set $C$ is said to be \emph{weakly positively invariant with respect
  to $T$} if for any $P_0 \in C$ there exists $P \in T(P_0)$ such that
$P \in C$.  A point $P_*$ is said to be a \emph{fixed point of $T$} if
$P_* \in T(P_*)$.  Let $U:Q^n \rightarrow \real$.  We say that $U$ is
a \emph{Lyapunov function for $T$ on $W$} if (i) $U$ is continuous on
$W$ and (ii) $U(P') \le U(P)$ for all $P' \in T(P)$ and all $P \in W$.

We now turn to the design of discrete-time algorithms for
limited-range coverage control. We start by extending the definition
of the aggregate objective function~$\HH$ to consider general partitions
$\WW$ of $Q$ as follows. Let $P \in Q^n$ and let $\WW= \left\{
  W_i\subset Q \right\}_{i =1}^n$ be a partition of $Q$ such that
$W_i$ is a convex polygon and $p_i \in W_i$, for $i \in \until{n}$.
Define the function
\begin{equation*}
  \HH_e(P, \WW) = \sum_{i=1}^n \int_{W_i} f(\| q - p_i \|) \phi(q) dq
  \, .  
\end{equation*}
The function $\HH_e$ is differentiable with respect to its first
variable for all $P \in Q^n$, and its partial derivative is given by
\begin{align}
  \pder{\HH_e}{p_i}(P,\WW) =&\; \sum_{\alpha=1}^{m+1} \int_{W_i \cap
    D_{(R_{\alpha-1},R_\alpha)}(p_i)} \pder{}{p_i} f_\alpha (\| q -
  p_i \| ) \phi(q) dq \nonumber \\
  &\;+ \sum_{\alpha=1}^{m+1} \int_{\partial (W_i \cap
    B_{R_\alpha}(p_i))} n_{i,\alpha}^t \pder{\gamma_{i,\alpha}}{p_i}
  f_\alpha (\dist(\cdot,p_i)) \, \phi \nonumber \\
  &\; - \sum_{\alpha=1}^{m+1} \int_{\partial (W_i \cap
    B_{R_{\alpha-1}}(p_i))} n_{i,\alpha-1}^t
  \pder{\gamma_{i,\alpha-1}}{p_i} f_\alpha (\dist(\cdot,p_i)) \, \phi
  \, ,
  \label{eq:gradient-fixed-partition}
\end{align} 
where for each $k \in \until{n}$ and each $\alpha \in \until{m+1}$,
$n_{k,\alpha}(q)$ denotes the unit outward normal to $W_k \cap
B_{R_\alpha}(p_k)$ at $q$, and $\gamma_{k,\alpha}:\sphere^1 \times Q^n
\rightarrow \real^2$ denotes a parameterization for the piecewise
smooth family $\left\{W_k \cap B_{R_\alpha}(p_k) \; | \; P \in Q^n
\right\}$. Note that, using the definition of $\HH_1$
(cf.~Section~\ref{sec:locational-optimization}), one can also write
\begin{align*}
  \HH_e(P, \WW) = \sum_{i=1}^n \HH_1(p_i,W_i) \, .
\end{align*}
The following two equalities hold
\begin{align}
  \HH(P) &= \HH_e (P,\VV (P)) \, , \quad \text{for all} \; P \in Q^n
  \, ,\\
  \pder{\HH_e}{p_i}(P,\VV(P)) & = \pder{\HH}{p_i} (P) \, , \quad
  \text{for all} \; P \in Q^n \setminus \Sc \, .
\end{align}

Let $P \in \Sc$ and consider a partition $\WW= \left\{ W_i\subset Q
\right\}_{i =1}^n$ of $Q$ such that $W_i$ is a convex polygon and $p_i
\in W_i$, for $i \in \until{n}$. Let $i_0,j_0 \in \until{n}$, $i_0
\neq j_0$ such that $p_{i_0} = p_{j_0}$. Then, following
Remark~\ref{re:Voronoi-covering}, $V_{i_0} (P) = V_{j_0} (P)$, and
$\VV (P)$ is no longer a partition of $Q$, but a covering.
Nevertheless, one could consider the line determined by the edge
$W_{i_0} \cap W_{j_0}$ and the associated hyperplanes $H_{i_0,W_{i_0}
  \cap W_{j_0}}$ and $H_{j_0,W_{i_0} \cap W_{j_0}}$ such that $W_{i_0}
\subset H_{i_0,W_{i_0} \cap W_{j_0}}$ and $W_{j_0} \subset
H_{j_0,W_{i_0} \cap W_{j_0}}$.  With a slight abuse of notation,
redefining
\begin{align*}
  V_{i_0}(P) = V_{i_0}(P) \cap H_{i_0,W_{i_0} \cap W_{j_0}} \, , \quad
  V_{j_0}(P) = V_{j_0}(P) \cap H_{j_0,W_{i_0} \cap W_{j_0}} \, ,
\end{align*}
the collection $\VV (P)$ can be seen a partition of~$Q$. This
procedure can be extended if there are more than two coincident agents
$\{i_1,\dots,i_s\}$ at a point $p \in Q$ by defining
\begin{align*}
  V_{i_\mu}(P) = V_{i_\mu}(P) \intersection \big(
  \intersection\limits_{\nu \in \until{s} \setminus\{ \mu\}}
  H_{i_\mu,W_{i_\mu} \cap W_{i_\nu}} \big) \, , \quad \mu
  \in \until{s} \, .
\end{align*}
In the following, such a construction will be tacitly performed
whenever we have a configuration $P \in \Sc$ and a partition $\WW$ of~$Q$.

The following lemma shows that the Voronoi partition is optimal within
the set of partitions of $Q$.
\begin{lmm}\label{le:fixed-partition}
  Let $\phi$ and $f$ be a density and a performance function,
  respectively.  Let $P \in Q^n$ and consider a partition $\WW=
  \left\{ W_i\subset Q \right\}_{i =1}^n$ of $Q$ such that $W_i$ is a
  convex polygon and $p_i \in W_i$, for $i \in \until{n}$. Then
  \begin{align*}
    \HH_e(P,\WW) \leq \HH_e(P,\VV(P)) \, ,
  \end{align*}
  and the inequality is strict if $f$ is strictly decreasing and the
  partitions $\VV(P)$ and $\WW$ differ by a set of non-zero measure.
\end{lmm}
\begin{proof}
  Given the chain of implications $q\in V_j(P) \Rightarrow
  \|q-p_i\|\geq\|q-p_j\| \Rightarrow f(\|q-p_i\|)\phi(q) \leq
  f(\|q-p_j\|)\phi(q)$, we compute
  \begin{align*}
    \HH_e(P,\WW) &=\; \sum_{i=1}^n\sum_{j=1}^n \int_{W_i\intersect
      V_j(P)} f(\|q-p_i\|)\phi(q)dq \\
    &\leq\; \sum_{i=1}^n\sum_{j=1}^n \int_{W_i\intersect V_j(P)}
    f(\|q-p_j\|)\phi(q)dq = \HH_e(P,\VV(P)).
  \end{align*}
\end{proof}

We are now ready to characterize a class of algorithms with guaranteed
convergence to the set of critical points of the aggregate objective function
$\HH$.
\begin{prpstn}[Discrete-time ascent]
  \label{prop:discrete-ascent}
  Let $T: Q^n \rightarrow 2^{Q^n}$ be an algorithm closed on $Q^n
  \setminus \Sc$ satisfying the following properties:
  \begin{itemize}
  \item[(a)] for all $P \in Q^n$, $T(P) \cap \Sc = \emptyset$;
  \item[(b)] for all $P \in Q^n \setminus \Sc$, $P' \in T(P)$ and $i
    \in \until{n}$, $\HH_1(p'_i,V_i(P)) \ge \HH_1(p_i,V_i(P))$;
  \item[(c)] for all $P \in \Sc$ and $P' \in T(P)$, $\HH(P') >
    \HH(P)$;
  \item[(d)] if $P \in Q^n \setminus \Sc$ is not a critical point of
    $\HH$, then for all $P' \in T(P)$, there exists $j \in \until{n}$
    such that $\HH_1(p'_j,V_j(P)) > \HH_1(p_j,V_j(P))$.
  \end{itemize}
  Let $P_0 \in Q^n$ denote the initial agents' location.  Then, any
  sequence $\setdef{P_n}{n\in\natural \cup \{0\}}$ generated according
  to equation~\eqref{eq:algorithm} converges to the set of critical
  points of~$\HH$.
\end{prpstn}

\begin{proof}
  Consider $-\HH: Q^n \rightarrow \real$ as a candidate Lyapunov
  function for the algorithm~$T$ on $Q^n \setminus \Sc$.  Because of
  Lemma~\ref{le:fixed-partition}, we have
  \[
  \HH (P') = \HH_e (P', \VV (P')) \ge \HH_e (P', \VV (P)) \, ,
  \]
  for all $P' \in T(P)$. In addition, because of property (b) of $T$, we
  also have
  \[
  \HH_e (P', \VV (P)) \ge \HH_e (P, \VV (P)) = \HH (P) \, ,
  \]
  for all $P' \in T(P)$.  Hence, $\HH (P') \le \HH (P)$ for all $P' \in
  T(P)$ and all $P \in Q^n$. Therefore, we deduce that $-\HH$ is a Lyapunov
  function for the algorithm $T$.  Let $P_0 \in Q^n \setminus \Sc$ and
  consider a sequence $\setdef{P_n}{n\in\natural \cup \{0\}}$ generated
  according to equation~\eqref{eq:algorithm}.  Because of property (a)
  of~$T$, $\setdef{P_n}{n\in\natural \cup \{0\}}$ remains in $Q^n \setminus
  \Sc \subset Q^n$. Since $Q^n$ is compact, we conclude that the sequence
  is bounded. Now, the application of the discrete-time LaSalle Invariance
  Principle (see Appendix~\ref{app:discrete-LaSalle},
  Theorem~\ref{th:discrete-LaSalle}) guarantees that there exists $c \in
  \real$ such that $P_n \rightarrow M \cap \HH^{-1}(c)$, where $M$ is the
  largest weakly positively invariant set contained in $\setdef{P' \in
    Q^n}{\exists P' \in T(P) \; \text{such that} \; \HH(P') = \HH (P)}$.
  Properties (c) and (d) of $T$ imply that $M$ must be contained in the set
  of critical points of $\HH$. If $P_0 \in \Sc$, the sequence
  $\setdef{P_n}{n\in\natural \cup \{0\}}$ can be equivalently described by
  $\{P_0\} \union \setdef{P_n}{n\in\natural}$.  Since $P_1 \in Q^n
  \setminus \Sc$ by property (a) of $T$, the previous argument implies that
  the sequence converges to the set of critical points of~$\HH$.
\end{proof}

In what follows, we devise a general algorithm $T:Q^n \rightarrow 2^{Q^n}$
verifying properties (a)-(d) in Proposition~\ref{prop:discrete-ascent}. We
shall do so by designing a discrete-time version of the gradient ascent
algorithm for continuous-time settings.

Recall that Lemma~\ref{le:maxima-not-in-boundary} asserts that if $p_0 \in
\partial W$, then $\pder{\HH_1(\cdot,W)}{p}(p_0) \neq 0$ points toward the
interior of $W$. If $p_0 \in \interior (W)$ is not a critical point, then
one also has that $\pder{\HH_1(\cdot,W)}{p}(p_0) \neq 0$.  For both cases,
there exists $\eps=\eps (p_0,W) >0$ such that the point $p_\delta$ defined
by
\begin{align*}
  p_\delta = p_0 + \delta \, \pder{\HH_1(\cdot,W)}{p}(p_0) \in W
\end{align*}
has the property that $\HH_1 (p_\delta) > \HH_1 (p_0)$, for all
$\delta \in (0,\eps)$, and $\HH_1 (p_\eps) = \HH_1 (p_0)$. As it is
usually done in nonlinear programming~\cite{DGL:84}, the computation
of the step-size $\eps$ can be implemented numerically via a ``line
search''.  With this discussion in mind, let us define the \emph{line
  search algorithm} $\subscr{T}{ls}:Q^n \rightarrow 2^{Q^n}$ as
follows:

\begin{quote}  
  Given $P \in Q^n$, let $P' \in \subscr{T}{ls}(P)$ if, for $i \in
  \until{n}$ with the property that $p_i \neq p_j, j \in \until{n}
  \setminus \{i\}$,
  \begin{align}\label{eq:line-search-algorithm-I}
    p'_i = p_i + \delta \, \pder{\HH_1(\cdot,V_i(P))}{p}(p_i) \, , \;
    \text{with} \; \delta \in \left[
      \frac{\eps(p_i,V_i(P))}{3},\frac{\eps(p_i,V_i(P))}{2} \right] ,
  \end{align}
  and, for each set $\{i_1, \dots, i_s \}$ of coincident indexes at a point
  $p \in Q$,
  \begin{align}\label{eq:line-search-algorithm-II}
    p'_{i_\mu} = p_{i_\mu} + \delta \,
    \pder{\HH_1(\cdot,Y_{i_\mu})}{p}(p_{i_\mu}) \, , \; \text{with}
    \; \delta \in \left[ \frac{\eps(p_i,Y_i)}{3},\frac{\eps(p_i,Y_i)}{2}
    \right] ,
  \end{align}
  where $\{ Y_{i_1},\dots,Y_{i_s} \}$ is a partition of $V_{i_1}(P) = \dots
  = V_{i_s}(P)$ verifying $p \in Y_{i_\mu}$, for $\mu \in \until{s}$.
\end{quote}

\begin{prpstn}\label{prop:line-search-algorithm}
  The algorithm $\subscr{T}{ls}:Q^n \rightarrow 2^{Q^n}$ defined by
  equations~\eqref{eq:line-search-algorithm-I}-\eqref{eq:line-search-algorithm-II}
  is closed on $Q^n \setminus \Sc$, and verifies properties (a)-(d) in
  Proposition~\ref{prop:discrete-ascent}.
\end{prpstn}

\begin{proof}
  The fact that $\subscr{T}{ls}$ is closed on $Q^n \setminus \Sc$
  follows from its definition and the continuous dependence of $\eps
  (p,V(P))$ on $P \in Q^n \setminus \Sc$. Regarding the properties in
  Proposition~\ref{prop:discrete-ascent}, consider the following
  discussion. Let $P \in Q^n$ and consider $P' \in \subscr{T}{ls}(P)$.
  On the one hand, equation~\eqref{eq:line-search-algorithm-I} and the
  definition of $\eps (p,V(P))$ implies that $p_i' \in \interior
  (V_i(P))$ for each $i \in \until{n}$ such that $p_i \neq p_j$ for
  all $j \in \until{n} \setminus \{i\}$. On the other hand,
  equation~\eqref{eq:line-search-algorithm-II} and
  Lemma~\ref{le:maxima-not-in-boundary} implies $p'_{i_\mu} \in
  \interior (Y_{i_\mu})$.  Therefore, we deduce that $P' \not \in
  \Sc$, and property (a) is verified.  Using
  equation~\eqref{eq:line-search-algorithm-I}, one has that for all $P
  \in Q^n \setminus \Sc$, $P' \in \subscr{T}{ls}(P)$ and all $i\in
  \until{n}$, $\HH_1(p'_i,V_i(P)) \ge \HH_1(p_i,V_i(P))$, i.e., the
  algorithm $\subscr{T}{ls}$ verifies property (b). With respect to
  property (c), let $P \in \Sc$. For simplicity, we only deal with the
  case when there exists $i,j \in \until{n}$, $i \neq j$ such that
  $p_i = p_j$, and all other $p_k \neq p_i=p_j$, $k \in \until{n}
  \setminus \{i,j\}$ are distinct among them (the cases with more
  degeneracies are treated analogously). Let $P' \in T(P)$.  According
  to equation~\eqref{eq:line-search-algorithm-II}, we have
  \begin{align*}
    \HH (P) = \sum_{k \in \until{n} \setminus \{i,j\}} \hspace*{-10pt}
    \HH_1(p_k,V_k(P)) + \HH_1(p_i,Y_i) + \HH_1(p_j,Y_j) \, ,
  \end{align*}
  where $\{ Y_{i},Y_j \}$ is a partition of $V_{i}(P) = V_{j}(P)$ with
  $p_i \in Y_{i}$ and $p_j \in Y_{j}$.  Since necessarily $p_i \in
  \partial Y_i$ and $p_j \in \partial Y_j$,
  Lemma~\ref{le:maxima-not-in-boundary} implies that $\HH_1(p_i,Y_i) +
  \HH_1(p_j,Y_j) < \HH_1(p'_i,Y_i) + \HH_1(p'_j,Y_j)$.  Therefore,
  $\HH (P) < \HH (P')$, i.e., property (c) is verified by
  $\subscr{T}{ls}$.  Finally, if $P \in Q^n \setminus \Sc$ is not a
  critical point of $\HH$, then there must exist $i \in \until{n}$
  such that
  \begin{align*}
    \pder{\HH}{p_i} (P) = \pder{\HH_e}{p_i}(P,\VV(P)) \neq 0 \, .
  \end{align*}
  Equivalently, $p_i$ is not a critical point of $\HH_1(\cdot,V_i(P))
  : V_i(P) \rightarrow \real$, and therefore $\eps(p_i,V_i(P)) > 0$.
  By equation~\eqref{eq:line-search-algorithm-I}, we conclude that
  $\HH_1(p'_i,V_i(P)) > \HH_1(p_i,V_i(P))$ for all $P' \in
  \subscr{T}{ls}(P)$, i.e., the algorithm $\subscr{T}{ls}$ verifies
  property (d).
\end{proof}

\begin{crllr}
  Consider the algorithm $\subscr{T}{ls}:Q^n \rightarrow 2^{Q^n}$
  defined by
  equations~\eqref{eq:line-search-algorithm-I}-\eqref{eq:line-search-algorithm-II}.
  Then
  \begin{enumerate}
  \item For a general performance function $f$, the algorithm
    $\subscr{T}{ls}$ is spatially distributed over the Delaunay
    graph~$\subscr{\GG}{D}(\PP)$.  If, in addition, the performance
    function verifies $f(x) = b$ for all $x \ge R$, then
    $\subscr{T}{ls}$ is spatially distributed over the $2R$-limited
    Delaunay graph $\subscr{\GG}{LD}(P,2R)$;
  \item The sequence of agents' locations generated by $\subscr{T}{ls}$
    according to equation~\eqref{eq:algorithm} starting at $P_0 \in
    Q^n$, converges asymptotically to the set of critical points of
    the aggregate objective function~$\HH$.
  \end{enumerate}
\end{crllr}
\begin{proof}
  The statement (i) is a direct consequence of
  Corollary~\ref{co:Delaunay-distributed}.  The convergence result is
  a consequence of Propositions~\ref{prop:discrete-ascent}
  and~\ref{prop:line-search-algorithm}.
\end{proof}

\begin{rmrk}
  As we noticed in Section~\ref{se:continuous-time}, in a practical
  setting, the network agents have typically a limited
  sensing/communication radius $r>0$, and therefore, following the
  result in Proposition~\ref{prop:approximation-HH}, we seek to
  maximize the function $\HH_{\frac{r}{2}}$. \oprocend
\end{rmrk}
    
In certain cases, it might be possible to construct specific algorithms
tailored to the concrete aggregate objective function at hand.  A relevant
example of this situation is when the local maxima of the function
$\HH_1(\cdot,W)$ can be characterized for each fixed polygon $W$.  With
this discussion in mind, let us define the \emph{max algorithm}
$\subscr{T}{max}:Q^n \rightarrow 2^{Q^n}$ as follows:
\begin{quote}
  For $P \in Q^n \setminus \Sc$, let
  \begin{align}\label{eq:maxima-algorithm}
    \subscr{T}{max}(P) & = \setdef{P' \in Q^n}{p'_i \; \text{is a local
        maximum of} \; \HH_1(\cdot,V_i(P)) \, , \; \text{for} \; i \in
      \until{n}} .
  \end{align}
  If $P \in \Sc$, for each set $\{i_1, \dots, i_s \}$ of coincident indexes
  at a point $p \in Q$, let $p'_{i_\mu}$ be a local maximum of
  $\HH_1(\cdot,Y_{i_\mu})$, where $\{ Y_{i_1},\dots,Y_{i_s} \}$ is a
  partition of $V_{i_1}(P) = \dots = V_{i_s}(P)$ verifying $p \in
  Y_{i_\mu}$, for $\mu \in \until{s}$.
\end{quote}
One can show that $\subscr{T}{max}$ is closed on $Q^n \setminus \Sc$
and verifies properties (a)-(d) in
Proposition~\ref{prop:discrete-ascent}.  As before, the algorithm
$\subscr{T}{max}$ is spatially distributed over the Delaunay
graph~$\subscr{\GG}{D}(\PP)$ and, if the performance function is
$f_{\frac{r}{2}}$, then $\subscr{T}{max}$ is spatially distributed
over the $r$-limited Delaunay graph~$\subscr{\GG}{LD}(\PP,r)$.
    
It is worth noticing that Lemma~\ref{le:maxima-not-in-boundary} guarantees
that the local maxima of $\HH_1(\cdot,W)$ are not in the boundary of $W$,
and therefore are contained in the set $\setdef{p_*\in
  W}{\pder{\HH_1(\cdot,W)}{p}(p_*)=0}$. Moreover, if $f$ is concave, then
$\HH_1$ is also concave, as stated in the following lemma.
\begin{lmm}
  If $\map{f}{\ov{\real}_+}{\real}$ is concave, then $\HH_1$ is
  concave.
\end{lmm}
\begin{proof}
  For fixed $q\in Q$, the map $p\mapsto f(\|q-p\|)\phi(q)$ is concave; the
  integral with respect to $q$ of a map with this property is concave in
  $p$; see~\cite[Subsection 3.2.1]{SB-LV:04}.
\end{proof}
As a consequence, the set of global maxima of $\HH_1 (\cdot,W)$ is compact,
convex and characterized by the equation
\begin{align*}
  \pder{\HH_1(\cdot,W)}{p}(p)=0 \, .
\end{align*}
In particular, these conditions are met in the centroid problem introduced
in Section~\ref{sec:locational-optimization}, where $f(x) = -x^2$ is
concave and the unique global minimum of $\HH_1(\cdot,W)$ is the centroid
$\CM_W$ of $W$. In this case, the algorithm $\subscr{T}{max}$ is precisely
the Lloyd quantization
algorithm~\cite{RMG-DLN:98,QD-VF-MG:99,JC-SM-TK-FB:02j}. 

\section{Simulations}\label{sec:simulations}

To illustrate the performance of the coordination algorithms, we
include some simulation results\footnote{Due to the limited size of
  the submissions to the arXiv e-print server, we could not include
  here the figures. The interested reader is referred to
  http:/\!/motion.csl.uiuc.edu for the complete preprint version with
  all figures.}. The algorithms are implemented in
\texttt{Mathematica} as a library of routines and a main program
running the simulation.  The objective of a first routine is to
compute the intersection of the bounded Voronoi cell $V_i$ with the
ball $B_{\frac{r}{2}}(p_i)$, for $i \in \until{n}$, and to
parameterize each set $V_i \cap B_{\frac{r}{2}}(p_i)$ in polar
coordinates. The objective of a second routine is to compute the
surface integrals on these sets and the line integrals on their
boundaries
via the numerical integration routine~\texttt{NIntegrate}.  We paid
careful attention to numerical accuracy issues in the computation of
the Voronoi diagram and in the integration.

We show executions of the discrete-time algorithm $\subscr{T}{ls}$
(cf.
equations~\eqref{eq:line-search-algorithm-I}-\eqref{eq:line-search-algorithm-II})
for the centroid problem, the area problem, the mixed centroid-area
problem with continuous sensing performance, and the mixed
centroid-area problem with discontinuous sensing performance in
Figs.~\ref{fig:coverage-centroid},~\ref{fig:coverage-area},~\ref{fig:coverage-mixed-cont}
and~\ref{fig:coverage-mixed-discont}, respectively. Measuring
displacements in meters, we consider the domain~$Q$ determined by the
vertexes
\begin{align*}
  \{ (0,0), (2.125, 0), (2.9325, 1.5), (2.975, 1.6), (2.9325, 1.7),
  (2.295, 2.1), (0.85, 2.3), (0.17, 1.2) \}.
\end{align*}
The diameter of this domain is $\diam{Q} = 3.37796$. In all examples,
the distribution density function~$\phi$ is the sum of five Gaussian
functions of the form $5 \, \exp(6(-(x-x_{\text{center}})^2 -
(y-y_{\text{center}})^2))$ and is represented by means of its contour
plot. The centers $(x_{\text{center}},y_{\text{center}})$ of the
Gaussians are given, respectively, by $(2,.25)$, $(1,2.25)$,
$(1.9,1.9)$, $(2.35,1.25)$ and $(.1,.1)$. Measured with respect
to~$\phi$, the area of the domain is~$\area{Q} = 8.61656$.  Unless
otherwise noted, each agent operates with a finite
sensing/communication radius equal to $r=.45$.

\begin{figure}[htbp!]
  \centering%
  \caption{Centroid problem (with agent
    performance $f(x) = -x^2$): discrete-time algorithm $\subscr{T}{ls}$
    for $16$ agents on a convex polygonal environment.  The left
    (respectively, right) figure illustrates the initial (respectively,
    final) locations and Voronoi partition.  The central figure illustrates
    the gradient ascent flow. After $90$~seconds, the value of the
    multi-center function is approximately~$-.321531$.}
  \label{fig:coverage-centroid}
\end{figure}
\begin{figure}[htbp!]
  \centering%
  \caption{Area problem (with agent performance $f(x) =
    \indicator{[0,\frac{r}{2}]}(x)$): discrete-time algorithm
    $\subscr{T}{ls}$ for $16$ agents on a convex polygonal
    environment.  The left (respectively, right) figure illustrates
    the initial (respectively, final) locations and Voronoi partition.
    The central figure illustrates the gradient ascent flow. For each
    agent $i$, the intersection $V_i \cap B_{\frac{r}{2}}(p_i)$ is
    plotted in light gray. After $18$~seconds, the value of the
    multi-center function is approximately $6.28977$.}
  \label{fig:coverage-area}  
\end{figure}
\begin{figure}[htbp!]
  \centering%
  \caption{Mixed centroid-area  problem (with continuous agent
    performance $f(x) = - x^2 \, \indicator{[0,\frac{r}{2})}(x) -
    \frac{r^2}{4} \cdot \indicator{[\frac{r}{2},+\infty)}(x)$):
    discrete-time algorithm $\subscr{T}{ls}$ for $16$ agents on a
    convex polygonal environment.  The left (respectively, right)
    figure illustrates the initial (respectively, final) locations and
    Voronoi partition.  The central figure illustrates the gradient
    ascent flow.  For each agent $i$, the intersection $V_i \cap
    B_{\frac{r}{2}}(p_i)$ is plotted in light gray. After
    $90$~seconds, the value of the multi-center function is
    approxima\-tely~$-0.252534$.}
  \label{fig:coverage-mixed-cont}  
\end{figure}
\begin{figure}[htbp!]
  \centering%
  \caption{Mixed centroid-area problem (with discontinuous agent
    performance $f(x) = - x^2 \, \indicator{[0,\frac{r}{2})}(x) -
    \diam{Q}^2 \cdot \indicator{[\frac{r}{2},+\infty)}(x)$):
    discrete-time algorithm $\subscr{T}{ls}$ for $16$ agents on a
    convex polygonal environment.  The left (respectively, right)
    figure illustrates the initial (respectively, final) locations and
    Voronoi partition.  The central figure illustrates the gradient
    ascent flow.  For each agent $i$, the intersection $V_i \cap
    B_{\frac{r}{2}}(p_i)$ is plotted in light gray. After
    $13.5$~seconds, the value of the multi-center function is
    approxima\-tely~$-6.803$.}
  \label{fig:coverage-mixed-discont}
\end{figure}

The execution of the coordination algorithm in
Figure~\ref{fig:coverage-mixed-discont} (with radius $r=.45$, agent
performance $f_{\frac{r}{2}}(x) = - x^2 \, \indicator{[0,\frac{r}{2})}(x) -
\diam{Q}^2 \cdot \indicator{[\frac{r}{2},+\infty)}(x)$ and corresponding
multi-center function $\HH_{\frac{r}{2}}$) can be regarded as a
limited-range implementation of the gradient ascent of the multi-center
function~$\HH$ corresponding to the agent performance $f(x) = -x^2$
(cf.~Figure~\ref{fig:coverage-centroid}); this performance function does
not have any range limitation. According to
Proposition~\ref{prop:approximation-HH}, we compute
\begin{align*}
  \beta & = \frac{f(\frac{r}{2})}{f (\diam{Q})} = \frac{1}{4} \left(
    \frac{r}{\diam{Q}} \right)^2 \approx 0.004437 \, , \\
  \Pi (P_\text{final}) & = \big( f \big(\tfrac{r}{2} \big)- f
  (\diam{Q}) \big) \area{Q \setminus \cup_{i=1}^n
    B_{\frac{r}{2}}(p_i)} \approx 26.5156 \, ,
\end{align*}
where $P_\text{final}$ denotes the final configuration in
Figure~\ref{fig:coverage-mixed-discont}.  From the constant-factor
approximation~\eqref{eq:sandwich-multiplicative} and the additive
approximation~\eqref{eq:sandwich-additive}, the absolute error is
guaranteed to be less than or equal to $\min \{(\beta-1)
\HH_{\frac{r}{2}} (P_{\text{final}}), \Pi (P_\text{final}) \}
\approx 6.77282$.
In order to compare the performance of this execution with the performance
of the discrete-time algorithm in the unlimited-range case, i.e., for the
case of $f(x)=-x^2$ (cf.  Figure~\ref{fig:coverage-centroid}), we compute
the percentage error in the value of the multi-center function~$\HH$ at
their final configurations. This percentage error is approximately equal
to~$30.7\%$.

Figure~\ref{fig:coverage-mixed-discont-II} below shows another
execution of the discrete-time algorithm $\subscr{T}{ls}$ for the
mixed centroid-area problem with discontinuous sensing performance,
where now the sensing/communication radius is taken equal to~$r=.65$.
In this case, the percentage error with respect to the performance of
the discrete-time algorithm in the unlimited-range case is
approximately equal to~$23\%$. As expected, the percentage error of
the performance of the limited-range implementation improves with
higher values of the ratio $\frac{r}{\diam{Q}}$.

\begin{figure}[htbp!]
  \centering%
  \caption{Execution of the discrete-time algorithm
    $\subscr{T}{ls}$ in the same setting as in
    Figure~\ref{fig:coverage-mixed-discont}, but with a
    sensing/communication radius equal to~$r=.65$. After $13.5$~seconds,
    the value of the multi-center function is approxima\-tely~$-1.10561$.}
  \label{fig:coverage-mixed-discont-II}
\end{figure}

\section{Conclusions and future work}\label{sec:conclusions}
We have presented novel spatially-distributed algorithms for coordinated
motion of groups of agents in continuous and discrete time.  Avenues of
possible future research include (1) distributed implementation of
deterministic annealing techniques~\cite{KR:98} (methods which promise to
overcome local maxima), (2) visibility-based algorithms for coverage in
non-convex environments, and (3) distributed algorithms for other
cooperative behaviors and sensing tasks, e.g., detection, estimation, and
map-building.

\begin{acknowledgement}
  This material is based upon work supported in part by ARO Grant
  DAAD~190110716, ONR YIP Award N00014-03-1-0512, and NSF SENSORS Award
  IIS-0330008. Sonia Mart{\'\i}nez's work was supported in part by a
  Fulbright PostDoctoral Fellowship from the Spanish Ministery of Education
  and Culture.
\end{acknowledgement}

\appendix

\section{Proof of Proposition~\ref{prop:mass-conservation}}
\label{sec:mass-conservation} 

\begin{proof}
  Let $x_0 \in (a,b)$.  Using the fact that the map $\gamma$ is
  continuous in both its arguments and that $\Omega(x_0)$ is strictly
  star-shaped, one can show that there exist an interval around $x_0$
  of the form $\mathcal{I}_{x_0} = (x_0-\eps,x_0+ \eps)$, a smooth
  function $\map{u_{x_0}}{\sphere^1 \times \ov{\real}_+}{\real^2}$ and
  a function $\map{r_{x_0}}{\sphere^1 \times
    \mathcal{I}_{x_0}}{\ov{\real}_+}$ smooth in $x$ and piecewise
  smooth in $\theta$ such that for all $x \in \mathcal{I}_{x_0}$, one
  has $\Omega(x)= \union_{\theta \in \sphere^1}
  \setdef{u_{x_0}(\theta,s)}{0 \le s \le r_{x_0}(\theta,x)}$ and
  $u_{x_0}(\theta,r_{x_0}(\theta,x))= \gamma(\theta,x)$, for $\theta
  \in \sphere^1$. For simplicity, we denote by $r$ and $u$ the
  functions $r_{x_0}$ and $u_{x_0}$, respectively.
  By definition, the function in~\eqref{eq:function} is continuously
  differentiable at $x_0$ if the following limit exists
  \begin{equation*}
    \lim_{h\rightarrow 0} \frac{1}{h} \left( \int_{\Omega(x_0+h)}
    \phi(q,x_0+h) dq - \int_{\Omega(x_0)} \phi(q,x_0) dq \right)  ,
  \end{equation*}
  and depends continuously on $x_0$.  Now, we can rewrite the previous
  limit as
  \begin{multline*}
    \lim_{h\rightarrow 0} \frac{1}{h} \int_0^{2 \pi} \left(
      \int_{0}^{r(\theta,x_0+h)} \!\! \phi(u(\theta,s),x_0+h) \Big \|
      \pder{u}{\theta} \! \times \! \pder{u}{s} \Big \| ds -
      \int_{0}^{r(\theta,x_0)} \phi(u(\theta,s),x_0) \Big \| \pder{u}{\theta}
      \times \pder{u}{s} \Big \| ds \right) d
    \theta = \\
    \lim_{h\rightarrow 0} \frac{1}{h} \int_0^{2 \pi} \left(
      \int_{r(\theta,x_0)}^{r(\theta,x_0+h)} \!\! \phi(u(\theta,s),x_0+h) \Big \|
      \pder{u}{\theta} \! \times \! \pder{u}{s} \Big \| ds
    \right. \\
    \left.  + \int_{0}^{r(\theta,x_0)} \hspace*{-5pt} \left(
        \phi(u(\theta,s),x_0+h) -\phi(u(\theta,s),x_0) \right) \Big \|
      \pder{u}{\theta} \times \pder{u}{s} \Big \| ds \right) d \theta \, ,
  \end{multline*}
  where $\times$ denotes the vector product and for brevity we omit
  that the partial derivatives $\pder{u}{\theta}$ and $\pder{u}{s}$
  are evaluated at $(\theta,s)$ in the integrals.  Now, since
  \begin{align*}
    \lim_{h\rightarrow 0} \frac{1}{h} \left( \phi(u(\theta,s),x_0+h)
      -\phi(u(\theta,s),x_0) \Big \| \pder{u}{\theta} \times \pder{u}{s}
      \Big \| \right) = \pder{\phi}{x_0} (u(\theta,s),x_0) \Big \|
    \pder{u}{\theta} \times \pder{u}{s} \Big \| \,
  \end{align*}
  almost everywhere and because this last function is measurable, the
  Lebesgue Dominated Convergence Theorem~\cite{RGB:95} implies that
  \begin{multline}\label{eq:aux-1}
    \lim_{h\rightarrow 0} \frac{1}{h} \int_0^{2 \pi} \hspace*{-5pt}
    \int_{0}^{r(\theta,x_0)} \left( \phi(u(\theta,s),x_0+h)
      -\phi(u(\theta,s),x_0) \right) \Big \| \pder{u}{\theta}
    \times  \pder{u}{s} \Big \|  ds  d \theta = \\
    \int_{0}^{2\pi} \hspace*{-5pt} \int_0^{r(\theta,x_0)}
    \pder{\phi}{x} (u(\theta,s),x_0) \Big \| \pder{u}{\theta} \times
    \pder{u}{s} \Big \| ds d \theta = \int_{\Omega(x_0)}
    \pder{\phi}{x}(q,x_0) dq \, .
  \end{multline}
  On the other hand, using the continuity of $\phi$, one can deduce that
  \begin{multline*}
    \lim_{h\rightarrow 0} \frac{1}{h} \int_0^{2 \pi} \hspace*{-5pt}
    \int_{r(\theta,x_0)}^{r(\theta,x_0+h)} \!\! \!
    \phi(u(\theta,s),x_0+h) \Big \| \pder{u}{\theta}(\theta,s) \!
    \times \! \pder{u}{s}(\theta,s) \Big
    \| ds \, d\theta \\
    = \lim_{h\rightarrow 0} \frac{1}{h} \int_0^{2 \pi} \hspace*{-5pt}
    \int_{x_0}^{x_0+h} \!\! \! \phi(u(\theta,r(\theta,z)),x_0+h) \Big
    \| \pder{u}{\theta}(\theta,r(\theta,z)) \! \times \!
    \pder{u}{s}(\theta,r(\theta,z)) \Big \| \pder{r}{x}(\theta,z) \,
    dz \,  d\theta \\
    = \int_0^{2 \pi} \phi(u(\theta,r(\theta,x_0)),x_0) \Big \|
    \pder{u}{\theta} (\theta,r(\theta,x_0)) \!  \times \!  \pder{u}{s}
    (\theta,r(\theta,x_0)) \Big \| \pder{r}{x_0}(\theta,x_0) \,
    d\theta \, .
  \end{multline*}
  Since $\gamma (\theta,x) = u(\theta,r(\theta,x))$ for all $\theta \in
  \sphere^1$ and $x \in \mathcal{I}_{x_0}$, one has
  \begin{align*}
    \pder{\gamma}{\theta}(\theta,x_0) &= \pder{u}{\theta}
    (\theta,r(\theta,x_0)) + \pder{u}{s} (\theta,r(\theta,x_0))
    \pder{r}{\theta} (\theta,x_0) \, , \\
    \pder{\gamma}{x} (\theta,x_0) &= \pder{u}{s}
    (\theta,r(\theta,x_0)) \pder{r}{x} (\theta,x_0) \, .
 \end{align*}
 Let $\chi$ denote the angle formed by $\pder{\gamma}{\theta}
 (\theta,x_0)$ and $\pder{u}{s} (\theta,r(\theta,x_0))$.  Then
 (omitting the expression $(\theta,r(\theta,x))$ for brevity),
  \begin{align*}
    \Big \| \pder{u}{\theta} \!  \times \!  \pder{u}{s} \Big \| = \Big \|
    \left( \pder{u}{\theta} + \pder{u}{s} \pder{r}{\theta} \right) \!
    \times \!  \pder{u}{s} \Big \| = \Big \| \der{\gamma}{\theta}\Big \|
    \Big \| \pder{u}{s} \Big \| \sin \chi = \Big \| \pder{\gamma}{\theta}
    \Big \| n^t (\gamma) \pder{u}{s} \, ,
  \end{align*}
  where in the last inequality we have used the fact that, since
  $\gamma_{x_0}$ is a parameterization of $\partial \Omega (x_0)$,
  then $\sin \chi = \cos \psi$, where~$\psi$ is the angle formed
  by $n$, the outward normal to $\partial \Omega (x_0)$, and
  $\pder{u}{s}$.  Therefore, we finally arrive at
  \begin{multline}\label{eq:aux-2}
    \int_0^{2 \pi} \phi(\gamma(\theta),x_0) \Big \| \pder{u}{\theta}
    (\theta,r(\theta,x_0)) \! \times \!  \pder{u}{s}
    (\theta,r(\theta,x_0)) \Big \| \pder{r}{x} (\theta,x_0) d\theta \\
    = \int_0^{2 \pi} \phi(\gamma(\theta),x_0) \Big \|
    \pder{\gamma}{\theta} (\theta,x_0) \Big \| n^t(\gamma(\theta,x_0))
    \pder{\gamma}{x} (\theta,x_0) d\theta = \int_{\partial
      \Omega(x_0)} \phi(\gamma,x_0) n^t(\gamma) \pder{\gamma}{x}
    d\gamma \, .
  \end{multline}
  Given the hypothesis of Proposition~\ref{prop:mass-conservation},
  both terms in~\eqref{eq:aux-1} and~\eqref{eq:aux-2} have a
  continuous dependence on $x_0 \in (a,b)$, which concludes the proof.
\end{proof}

\section{Upper bound on the area of the intersection between two balls}
\label{app:overlap-area}
\begin{lmm}\label{le:overlap-area}
  For $R\in\real_+$, let $p$, $p' \in\real^2$ satisfy $\|p -p'\| \le
  R$.  Then the area $\text{A}$ of $B_R(p')\intersection B_R^c(p)$
  satisfies $\text{A} \le \tfrac{2\sqrt{3} + 3}{3} R \|p -p' \|$.
\end{lmm}
\begin{proof}
  The area $\text{A}$ equals $\pi \, R^2 - \text{L}$, where $\text{L}$
  is the area of the non-trivial lune $B_R(p)\intersection B_R(p')$
  (see Figure~\ref{fig:lune}). 
  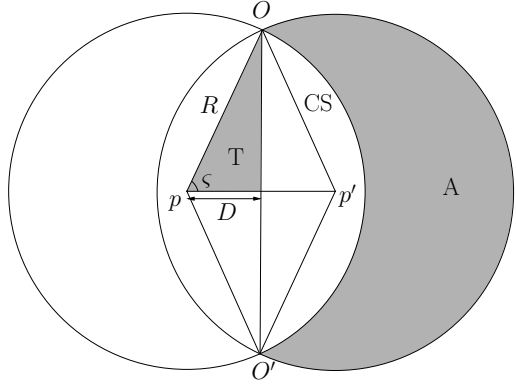
\begin{figure}[htb] 
    \begin{center}
      \resizebox{.4\linewidth}{!}{\input{area.tex}}
    \end{center}
    \caption{Areas of interest: $\text{A}$ is the area of
      $B_R(p')\intersection B_R^c(p)$, $\text{T}$ is the area of the
      triangle $T(p,O, \frac{p' +p}{2})$, and $\text{CS}$ is the area of
      the circular sector inside $B_R(p)$ determined by $(p,O,O')$.}
    \label{fig:lune}
  \end{figure}%
  Let $O$, $O'$ be the two points in the intersection $\partial B_R(p)
  \cap \partial B_R(p')$. Note that the triangle with vertices $O$,
  $p$ and $p'$, $T(O,p,p')$, and the triangle $T(O',p,p')$ are
  isosceles. This implies that the diagonals of the polygon
  $P(p,O,p',O')$ intersect at the middle point $\frac{p' + p}{2}$. Let
  $\varsigma$ be the angle of $T(O,p,p')$ at $p$ and $D= \frac{\| p -
    p'\|}{2}$. Then, the area of the lune $\text{L}$ can be computed
  as $\text{L} = 2(\text{CS} -2 \text{T})$, where $\text{CS}$ is the
  area of a circular sector with angle $2 \varsigma$ and $\text{T}$ is the
  area of the triangle $T(p,O, \frac{p' +p}{2})$.  Since $\text{CS} =
  \varsigma R^2$ and $\text{T}= \frac{1}{2}D \sqrt{R^2 -D^2}$, we have
  that $\text{A} = R^2 (\pi - 2 \varsigma) + 2 D \sqrt{R^2 -D^2}$. Now,
  using that $0 \le D \le R$, we deduce that
  \begin{equation*}
    \text{A} \le  R^2 (\pi - 2 \varsigma) + 2 R D\,.
  \end{equation*}
  In order to finally bound the first term of the right-hand side of
  the inequality with a quantity proportional to $D$, we use that
  $\varsigma = \arccos (\frac{D}{R})$. Consider now the function
  $g(x)= \pi -2 \arccos x - K x$. It is easy to see that for $K >
  \frac{4}{\sqrt{3}}$, one has $g'(x) \le 0$ and $g(x) \le 0$ when $0
  \le x \le \frac{1}{2}$. In particular this implies that, for $0 \le
  \frac{D}{R} \le \frac{1}{2}$, we have that $\pi - 2 \arccos
  (\frac{D}{R}) \le K \frac{D}{R}$. In other words, the former
  inequality is valid for $\| p -p'\|\le R$. This concludes the proof.
\end{proof}

\section{Discrete-time LaSalle Invariance
  Principle}\label{app:discrete-LaSalle}

The following result is an extension of two classical results: on the
one hand, it extends the discrete-time version of LaSalle Invariance
Principle~\cite{JPL:86} to algorithms defined via set-valued maps.  On
the other hand, it considers a more general notion of Lyapunov
function (cf.~Section~\ref{se:discrete-time}) than in the usual
statement of the Global Convergence Theorem~\cite{DGL:84}.

\begin{thrm}[Discrete-time LaSalle Invariance
  Principle]\label{th:discrete-LaSalle} Let $T$ be a closed algorithm
  on $W \subset \real^N$ and let $U$ be a Lyapunov function for $T$ on
  $W$. Let $x_0 \in W$ and assume the sequence $\setdef{x_n}{n \in
    \natural \cup \{0\}}$ defined via $x_{n+1} \in T(x_n)$ is in $W$
  and bounded.  Then there exists $c \in \real$ such that
  \begin{align*}
    x_n \longrightarrow M \cap U^{-1}(c) \, ,
  \end{align*}
  where $M$ is the largest weakly positively invariant set contained
  in 
  \begin{equation*}
  \setdef{x \in \real^N}{\exists y \in T(x) \; \text{such that} \;
    U(y)=U(x)} \cap \ov{W}.    
  \end{equation*}
\end{thrm}

\begin{proof}
  Let $\Omega (x_n) \subset \ov{W}$ denote the $\omega$-limit set of
  the sequence $\setdef{x_n}{n \in \natural \cup \{0\}}$. First, let
  us prove that $\Omega (x_n)$ is weakly positively invariant. Let $x
  \in \Omega (x_n)$. Then there exists a subsequence
  $\setdef{x_{n_m}}{m \in \natural \cup \{0\}}$ of $\setdef{x_n}{n \in
    \natural \cup \{0\}}$ such that $x_{n_m} \rightarrow x$.  Consider
  the sequence $\{x_{n_m + 1} \; | \; m \in \natural \cup \{0\}\}$.
  Since this sequence is bounded, it has a convergent subsequence. For
  ease of notation, we use the same notation to refer to it, i.e.,
  there exits $y$ such that $x_{n_m + 1} \rightarrow y$.  By
  definition, $y \in \Omega (x_n)$. Moreover, using the fact that $T$
  is closed, we deduce that $y \in T(x)$. Therefore $\Omega (x_n)$ is
  weakly positively invariant.
  
  Now, consider the sequence $\setdef{U(x_n)}{n \in \natural \cup \{0\}}$.
  Since $\setdef{x_n}{n \in \natural \cup \{0\}}$ is bounded and $U$ is a
  Lyapunov function for $T$ on $W$, this sequence is decreasing and bounded
  from below, and therefore convergent. Let $c \in \real$ such that $U(x_n)
  \rightarrow c$. Let us see that the value of $U$ on $\Omega (x_n)$ is
  constant and equal to $c$. Take any $x \in \Omega (x_n)$. Accordingly,
  there exists a subsequence $\setdef{x_{n_m}}{m \in \natural \cup \{0\}}$
  such that $x_{n_m} \rightarrow x$. Since $U$ is continuous, $U(x_{n_m})
  \rightarrow U(x)$. From $U(x_n) \rightarrow c$, we conclude that
  $U(x)=c$.
  
  Finally, the fact that $\Omega (x_n)$ is weakly positively invariant
  and $U$ is constant on $\Omega (x_n)$, implies that
  \[
  \Omega (x_n) \subset \setdef{x \in \real^N}{\exists y \in T(x) \;
    \text{such that} \; U(y) = U(x)} .
  \]
  Therefore, we conclude that $x_n \rightarrow M \cap U^{-1}(c)$,
  where $M$ is the largest weakly positively invariant set contained
  in $\setdef{x \in \real^N}{\exists y \in T(x) \; \text{such that} \;
    U(y) = U(x)} \cap \ov{W}$.
\end{proof}

\end{document}

%% file: counter-example.tex
\begin{picture}(0,0)%
\includegraphics{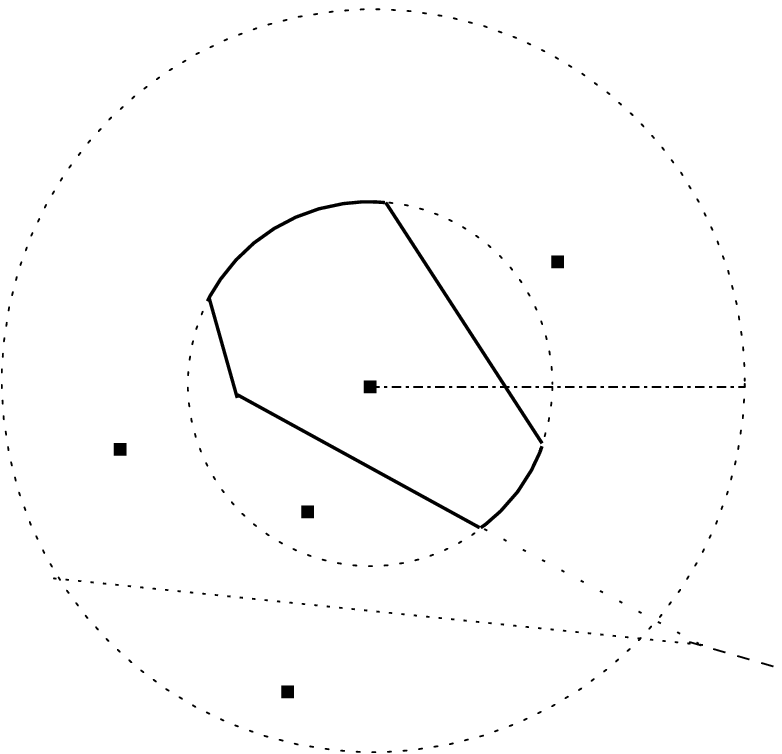}%
\end{picture}%
\setlength{\unitlength}{3947sp}%
\begingroup\makeatletter\ifx\SetFigFont\undefined%
\gdef\SetFigFont#1#2#3#4#5{%
  \reset@font\fontsize{#1}{#2pt}%
  \fontfamily{#3}\fontseries{#4}\fontshape{#5}%
  \selectfont}%
\fi\endgroup%
\begin{picture}(3732,3582)(3025,-5722)
\put(5791,-3496){\makebox(0,0)[lb]{\smash{\SetFigFont{17}{20.4}{\rmdefault}{\mddefault}{\updefault}{\color[rgb]{0,0,0}$p_j$}%
}}}
\put(4869,-4148){\makebox(0,0)[lb]{\smash{\SetFigFont{17}{20.4}{\rmdefault}{\mddefault}{\updefault}{\color[rgb]{0,0,0}$p_i$}%
}}}
\put(4606,-4717){\makebox(0,0)[lb]{\smash{\SetFigFont{17}{20.4}{\rmdefault}{\mddefault}{\updefault}{\color[rgb]{0,0,0}$p_k$}%
}}}
\put(5754,-4224){\makebox(0,0)[lb]{\smash{\SetFigFont{17}{20.4}{\rmdefault}{\mddefault}{\updefault}{\color[rgb]{0,0,0}$\frac{r}{2}$}%
}}}
\put(6677,-4163){\makebox(0,0)[lb]{\smash{\SetFigFont{17}{20.4}{\rmdefault}{\mddefault}{\updefault}{\color[rgb]{0,0,0}$r$}%
}}}
\put(3278,-4111){\makebox(0,0)[lb]{\smash{\SetFigFont{17}{20.4}{\rmdefault}{\mddefault}{\updefault}{\color[rgb]{0,0,0}$p_m$}%
}}}
\put(4562,-5571){\makebox(0,0)[lb]{\smash{\SetFigFont{17}{20.4}{\rmdefault}{\mddefault}{\updefault}{\color[rgb]{0,0,0}$p_l$}%
}}}
\end{picture}

%% file: restricted-Delaunay.tex
\begin{picture}(0,0)%
\includegraphics{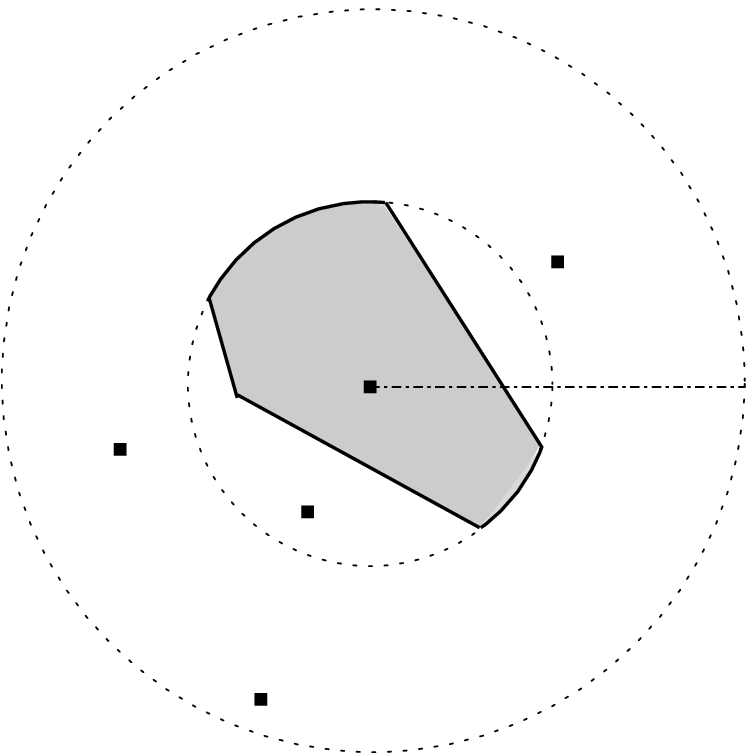}%
\end{picture}%
\setlength{\unitlength}{3947sp}%
\begingroup\makeatletter\ifx\SetFigFont\undefined%
\gdef\SetFigFont#1#2#3#4#5{%
  \reset@font\fontsize{#1}{#2pt}%
  \fontfamily{#3}\fontseries{#4}\fontshape{#5}%
  \selectfont}%
\fi\endgroup%
\begin{picture}(3652,3582)(3025,-5722)
\put(5791,-3496){\makebox(0,0)[lb]{\smash{\SetFigFont{17}{20.4}{\rmdefault}{\mddefault}{\updefault}{\color[rgb]{0,0,0}$p_j$}%
}}}
\put(4139,-3796){\makebox(0,0)[lb]{\smash{\SetFigFont{17}{20.4}{\rmdefault}{\mddefault}{\updefault}{\color[rgb]{0,0,0}$\Delta_{im}$}%
}}}
\put(5556,-4584){\makebox(0,0)[lb]{\smash{\SetFigFont{17}{20.4}{\rmdefault}{\mddefault}{\updefault}{\color[rgb]{0,0,0}$\operatorname{arc}_{i,1}$}%
}}}
\put(4869,-4148){\makebox(0,0)[lb]{\smash{\SetFigFont{17}{20.4}{\rmdefault}{\mddefault}{\updefault}{\color[rgb]{0,0,0}$p_i$}%
}}}
\put(4606,-4717){\makebox(0,0)[lb]{\smash{\SetFigFont{17}{20.4}{\rmdefault}{\mddefault}{\updefault}{\color[rgb]{0,0,0}$p_k$}%
}}}
\put(5754,-4224){\makebox(0,0)[lb]{\smash{\SetFigFont{17}{20.4}{\rmdefault}{\mddefault}{\updefault}{\color[rgb]{0,0,0}$\frac{r}{2}$}%
}}}
\put(6677,-4163){\makebox(0,0)[lb]{\smash{\SetFigFont{17}{20.4}{\rmdefault}{\mddefault}{\updefault}{\color[rgb]{0,0,0}$r$}%
}}}
\put(3884,-3002){\makebox(0,0)[lb]{\smash{\SetFigFont{17}{20.4}{\rmdefault}{\mddefault}{\updefault}{\color[rgb]{0,0,0}$\operatorname{arc}_{i,2}$}%
}}}
\put(4348,-4396){\makebox(0,0)[lb]{\smash{\SetFigFont{17}{20.4}{\rmdefault}{\mddefault}{\updefault}{\color[rgb]{0,0,0}$\Delta_{ik}$}%
}}}
\put(4790,-3631){\makebox(0,0)[lb]{\smash{\SetFigFont{17}{20.4}{\rmdefault}{\mddefault}{\updefault}{\color[rgb]{0,0,0}$\Delta_{ij}$}%
}}}
\put(4441,-5572){\makebox(0,0)[lb]{\smash{\SetFigFont{17}{20.4}{\rmdefault}{\mddefault}{\updefault}{\color[rgb]{0,0,0}$p_l$}%
}}}
\put(3278,-4111){\makebox(0,0)[lb]{\smash{\SetFigFont{17}{20.4}{\rmdefault}{\mddefault}{\updefault}{\color[rgb]{0,0,0}$p_m$}%
}}}
\end{picture}

%% file: restricted-Delaunay-2.tex
\begin{picture}(0,0)%
\includegraphics{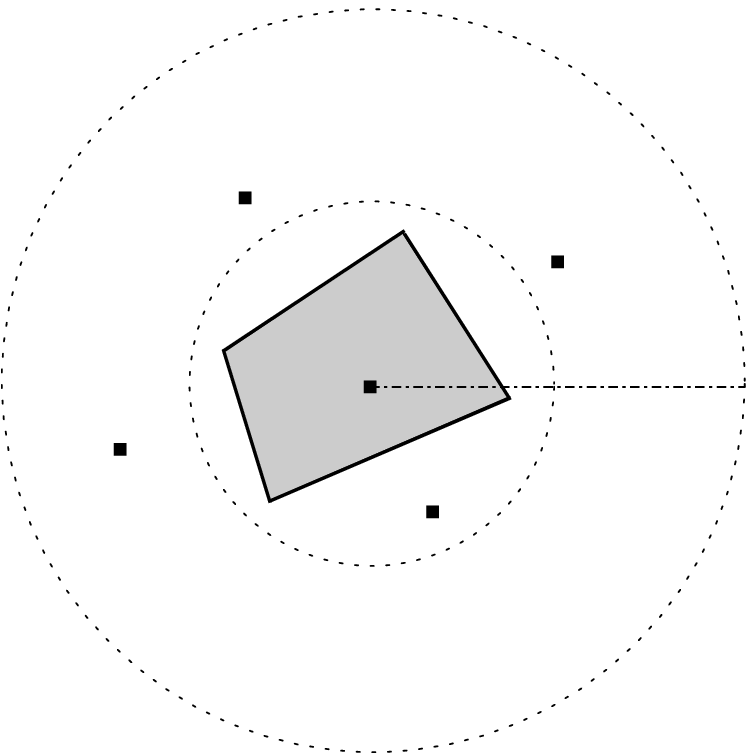}%
\end{picture}%
\setlength{\unitlength}{3947sp}%
\begingroup\makeatletter\ifx\SetFigFont\undefined%
\gdef\SetFigFont#1#2#3#4#5{%
  \reset@font\fontsize{#1}{#2pt}%
  \fontfamily{#3}\fontseries{#4}\fontshape{#5}%
  \selectfont}%
\fi\endgroup%
\begin{picture}(3652,3582)(3025,-5722)
\put(5791,-3496){\makebox(0,0)[lb]{\smash{\SetFigFont{17}{20.4}{\rmdefault}{\mddefault}{\updefault}{\color[rgb]{0,0,0}$p_j$}%
}}}
\put(4869,-4148){\makebox(0,0)[lb]{\smash{\SetFigFont{17}{20.4}{\rmdefault}{\mddefault}{\updefault}{\color[rgb]{0,0,0}$p_i$}%
}}}
\put(5754,-4224){\makebox(0,0)[lb]{\smash{\SetFigFont{17}{20.4}{\rmdefault}{\mddefault}{\updefault}{\color[rgb]{0,0,0}$\frac{r}{2}$}%
}}}
\put(6677,-4163){\makebox(0,0)[lb]{\smash{\SetFigFont{17}{20.4}{\rmdefault}{\mddefault}{\updefault}{\color[rgb]{0,0,0}$r$}%
}}}
\put(4790,-3631){\makebox(0,0)[lb]{\smash{\SetFigFont{17}{20.4}{\rmdefault}{\mddefault}{\updefault}{\color[rgb]{0,0,0}$\Delta_{ij}$}%
}}}
\put(3308,-4104){\makebox(0,0)[lb]{\smash{\SetFigFont{17}{20.4}{\rmdefault}{\mddefault}{\updefault}{\color[rgb]{0,0,0}$p_m$}%
}}}
\put(3916,-2902){\makebox(0,0)[lb]{\smash{\SetFigFont{17}{20.4}{\rmdefault}{\mddefault}{\updefault}{\color[rgb]{0,0,0}$p_l$}%
}}}
\put(4206,-4021){\makebox(0,0)[lb]{\smash{\SetFigFont{17}{20.4}{\rmdefault}{\mddefault}{\updefault}{\color[rgb]{0,0,0}$\Delta_{im}$}%
}}}
\put(4228,-3436){\makebox(0,0)[lb]{\smash{\SetFigFont{17}{20.4}{\rmdefault}{\mddefault}{\updefault}{\color[rgb]{0,0,0}$\Delta_{il}$}%
}}}
\put(4726,-4717){\makebox(0,0)[lb]{\smash{\SetFigFont{17}{20.4}{\rmdefault}{\mddefault}{\updefault}{\color[rgb]{0,0,0}$p_k$}%
}}}
\put(5045,-4389){\makebox(0,0)[lb]{\smash{\SetFigFont{17}{20.4}{\rmdefault}{\mddefault}{\updefault}{\color[rgb]{0,0,0}$\Delta_{ik}$}%
}}}
\end{picture}

%% file: area.tex
\begin{picture}(0,0)%
\includegraphics{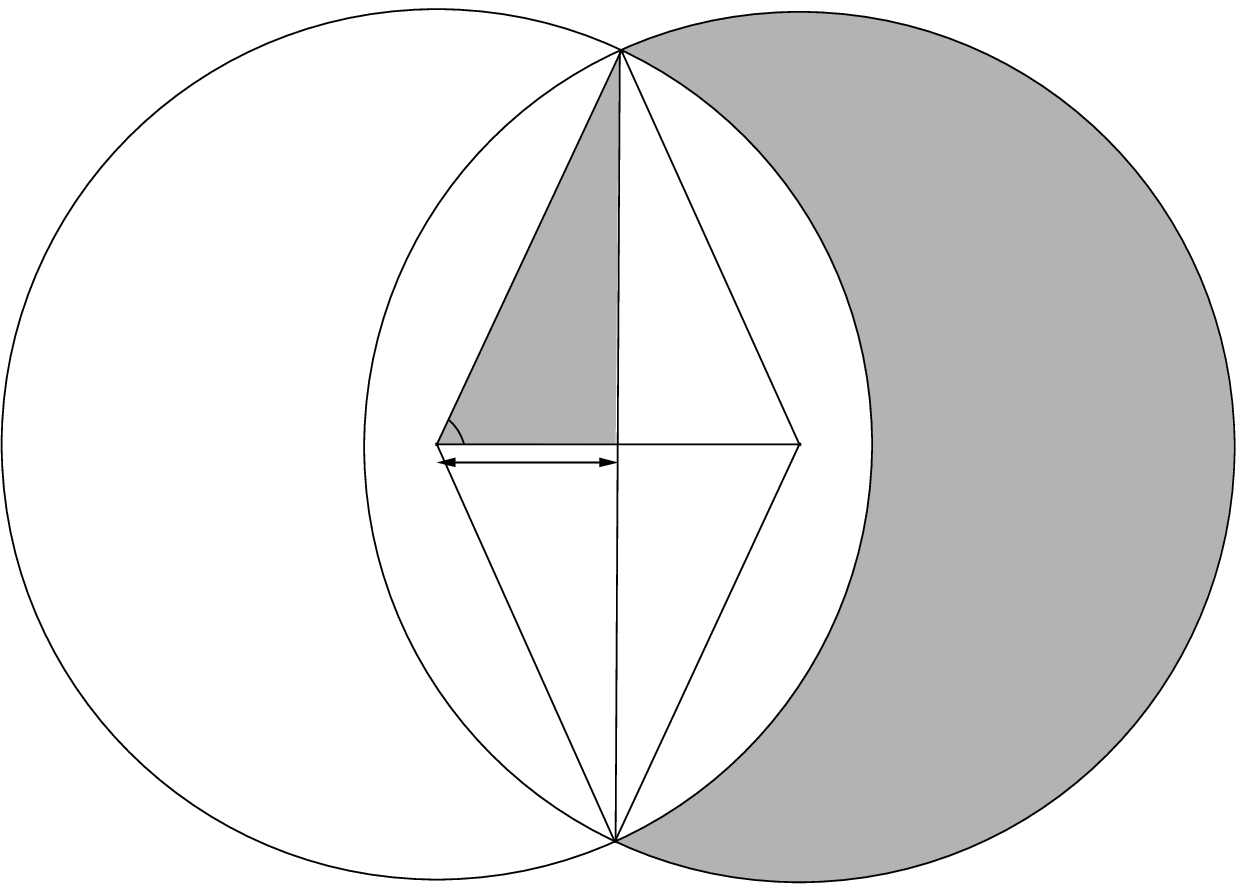}%
\end{picture}%
\setlength{\unitlength}{3947sp}%
\begingroup\makeatletter\ifx\SetFigFont\undefined%
\gdef\SetFigFont#1#2#3#4#5{%
  \reset@font\fontsize{#1}{#2pt}%
  \fontfamily{#3}\fontseries{#4}\fontshape{#5}%
  \selectfont}%
\fi\endgroup%
\begin{picture}(5934,4507)(2736,-5415)
\put(6616,-3267){\makebox(0,0)[lb]{\smash{\SetFigFont{17}{20.4}{\rmdefault}{\mddefault}{\updefault}{\color[rgb]{0,0,0}$p'$}%
}}}
\put(5311,-2832){\makebox(0,0)[lb]{\smash{\SetFigFont{17}{20.4}{\rmdefault}{\mddefault}{\updefault}{\color[rgb]{0,0,0}$\text{T}$}%
}}}
\put(7834,-3180){\makebox(0,0)[lb]{\smash{\SetFigFont{17}{20.4}{\rmdefault}{\mddefault}{\updefault}{\color[rgb]{0,0,0}$\text{A}$}%
}}}
\put(6202,-2197){\makebox(0,0)[lb]{\smash{\SetFigFont{17}{20.4}{\rmdefault}{\mddefault}{\updefault}{\color[rgb]{0,0,0}$\text{CS}$}%
}}}
\put(5007,-3063){\makebox(0,0)[lb]{\smash{\SetFigFont{17}{20.4}{\rmdefault}{\mddefault}{\updefault}{\color[rgb]{0,0,0}$\varsigma$}%
}}}
\put(5181,-3472){\makebox(0,0)[lb]{\smash{\SetFigFont{17}{20.4}{\rmdefault}{\mddefault}{\updefault}{\color[rgb]{0,0,0}$D$}%
}}}
\put(5591,-5343){\makebox(0,0)[lb]{\smash{\SetFigFont{17}{20.4}{\rmdefault}{\mddefault}{\updefault}{\color[rgb]{0,0,0}$O'$}%
}}}
\put(4983,-2215){\makebox(0,0)[lb]{\smash{\SetFigFont{17}{20.4}{\rmdefault}{\mddefault}{\updefault}{\color[rgb]{0,0,0}$R$}%
}}}
\put(5598,-1100){\makebox(0,0)[lb]{\smash{\SetFigFont{17}{20.4}{\rmdefault}{\mddefault}{\updefault}{\color[rgb]{0,0,0}$O$}%
}}}
\put(4626,-3310){\makebox(0,0)[lb]{\smash{\SetFigFont{17}{20.4}{\rmdefault}{\mddefault}{\updefault}{\color[rgb]{0,0,0}$p$}%
}}}
\end{picture}

%% file: main.bbl
\begin{thebibliography}{10}
\providecommand{\url}[1]{#1}
\csname url@rmstyle\endcsname
\providecommand{\newblock}{\relax}
\providecommand{\bibinfo}[2]{#2}
\providecommand\BIBentrySTDinterwordspacing{\spaceskip=0pt\relax}
\providecommand\BIBentryALTinterwordstretchfactor{4}
\providecommand\BIBentryALTinterwordspacing{\spaceskip=\fontdimen2\font plus
\BIBentryALTinterwordstretchfactor\fontdimen3\font minus
  \fontdimen4\font\relax}
\providecommand\BIBforeignlanguage[2]{{%
\expandafter\ifx\csname l@#1\endcsname\relax
\typeout{** WARNING: IEEEtran.bst: No hyphenation pattern has been}%
\typeout{** loaded for the language `#1'. Using the pattern for}%
\typeout{** the default language instead.}%
\else
\language=\csname l@#1\endcsname
\fi
#2}}

\bibitem{CWR:87}
C.~W. Reynolds, ``Flocks, herds, and schools: A distributed behavioral model,''
  \emph{Computer Graphics}, vol.~21, no.~4, pp. 25--34, 1987.

\bibitem{AO:86}
A.~Okubo, ``Dynamical aspects of animal grouping: swarms, schools, flocks and
  herds,'' \emph{Advances in Biophysics}, vol.~22, pp. 1--94, 1986.

\bibitem{RCA:98}
R.~C. Arkin, \emph{Behavior-Based Robotics}.\hskip 1em plus 0.5em minus
  0.4em\relax New York, NY: Cambridge University Press, 1998.

\bibitem{AJ-JL-ASM:02}
A.~Jadbabaie, J.~Lin, and A.~S. Morse, ``Coordination of groups of mobile
  autonomous agents using nearest neighbor rules,'' \emph{IEEE Transactions on
  Automatic Control}, vol.~48, no.~6, pp. 988--1001, 2003.

\bibitem{ROS-RMM:03c}
R.~Olfati-Saber and R.~M. Murray, ``Consensus problems in networks of agents
  with switching topology and time-delays,'' \emph{IEEE Transactions on
  Automatic Control}, Apr. 2003, submitted.

\bibitem{HT-AJ-GJP:03c}
H.~Tanner, A.~Jadbabaie, and G.~J. Pappas, ``Flocking in fixed and switching
  networks,'' \emph{IFAC Automatica}, July 2003, submitted.

\bibitem{PO-EF-NEL:03}
P.~\"Ogren, E.~Fiorelli, and N.~E. Leonard, ``Cooperative control of mobile
  sensor networks: adaptive gradient climbing in a distributed environment,''
  \emph{IEEE Transactions on Automatic Control}, July 2003, submitted.

\bibitem{KMP:04}
K.~M. Passino, \emph{Biomimicry for Optimization, Control, and
  Automation}.\hskip 1em plus 0.5em minus 0.4em\relax New York, NY: Springer
  Verlag, 2004, in print.

\bibitem{JC-SM-TK-FB:02j}
J.~Cort{\'e}s, S.~Mart{\'\i}nez, T.~Karatas, and F.~Bullo, ``Coverage control
  for mobile sensing networks,'' \emph{IEEE Transactions on Robotics and
  Automation}, 2003, to appear.

\bibitem{JC-FB:02m}
J.~Cort{\'e}s and F.~Bullo, ``Coordination and geometric optimization via
  distributed dynamical systems,'' \emph{SIAM Journal on Control and
  Optimization}, May 2003, submitted.

\bibitem{RD:00}
R.~Diestel, \emph{Graph Theory}, 2nd~ed., ser. Graduate Texts in
  Mathematics.\hskip 1em plus 0.5em minus 0.4em\relax New York, NY: Springer
  Verlag, 2000, vol. 173.

\bibitem{AO-BB-KS-SNC:00}
A.~Okabe, B.~Boots, K.~Sugihara, and S.~N. Chiu, \emph{Spatial Tessellations:
  Concepts and Applications of Voronoi Diagrams}, 2nd~ed., ser. Wiley Series in
  Probability and Statistics.\hskip 1em plus 0.5em minus 0.4em\relax New York,
  NY: John Wiley \& Sons, 2000.

\bibitem{ZD-HWH:01}
Z.~Drezner and H.~W. Hamacher, Eds., \emph{Facility Location: Applications and
  Theory}.\hskip 1em plus 0.5em minus 0.4em\relax New York, NY: Springer
  Verlag, 2001.

\bibitem{UH-JBM:94}
U.~Helmke and J.~Moore, \emph{Optimization and Dynamical Systems}.\hskip 1em
  plus 0.5em minus 0.4em\relax New York, NY: Springer Verlag, 1994.

\bibitem{JWJ-GTT:92}
J.~W. Jaromczyk and G.~T. Toussaint, ``Relative neighborhood graphs and their
  relatives,'' \emph{Proceedings of the IEEE}, vol.~80, no.~9, pp. 1502--1517,
  1992.

\bibitem{QD-VF-MG:99}
Q.~Du, V.~Faber, and M.~Gunzburger, ``Centroidal {V}oronoi tessellations:
  applications and algorithms,'' \emph{SIAM Review}, vol.~41, no.~4, pp.
  637--676, 1999.

\bibitem{AO-AS:97}
A.~Okabe and A.~Suzuki, ``Locational optimization problems solved through
  {V}oronoi diagrams,'' \emph{European Journal of Operational Research},
  vol.~98, no.~3, pp. 445--56, 1997.

\bibitem{RMG-DLN:98}
R.~M. Gray and D.~L. Neuhoff, ``Quantization,'' \emph{IEEE Transactions on
  Information Theory}, vol.~44, no.~6, pp. 2325--2383, 1998, {Commemorative
  Issue 1948-1998}.

\bibitem{YA:91}
Y.~Asami, ``A note on the derivation of the first and second derivative of
  objective functions in geographical optimization problems,'' \emph{Journal of
  the Faculty of Engineering, The University of Tokio (B)}, vol. XLI, no.~1,
  pp. 1--13, 1991.

\bibitem{MdB-MvK-MO:97}
M.~de~Berg, M.~van Kreveld, and M.~Overmars, \emph{Computational Geometry:
  Algorithms and Applications}.\hskip 1em plus 0.5em minus 0.4em\relax New
  York, NY: Springer Verlag, 1997.

\bibitem{XYL:03}
X.~Li, ``Algorithmic, geometric and graphs issues in wireless networks,''
  \emph{Wireless Communications and Mobile Computing}, vol.~3, no.~2, pp.
  119--140, 2003.

\bibitem{JG-LJG-JH-LZ-AZ:01}
J.~Gao, L.~J. Guibas, J.~Hershberger, L.~Zhang, and A.~Zhu, ``Geometric spanner
  for routing in mobile networks,'' in \emph{ACM International Symposium on
  Mobile Ad-hoc Networking \& Computing}, Long Beach, CA, Oct. 2001, pp.
  45--55.

\bibitem{AJC-JEM:94}
A.~J. Chorin and J.~E. Marsden, \emph{A Mathematical Introduction to Fluid
  Mechanics}, 3rd~ed., ser. Texts in Applied Mathematics.\hskip 1em plus 0.5em
  minus 0.4em\relax New York, NY: Springer Verlag, 1994, vol.~4.

\bibitem{HKK:96}
H.~K. Khalil, \emph{Nonlinear Systems}, 2nd~ed.\hskip 1em plus 0.5em minus
  0.4em\relax Englewood Cliffs, NJ: Prentice Hall, 1995.

\bibitem{DGL:84}
D.~G. Luenberger, \emph{Linear and Nonlinear Programming}, 2nd~ed.\hskip 1em
  plus 0.5em minus 0.4em\relax Reading, MA: Addison-Wesley, 1984.

\bibitem{SB-LV:04}
S.~Boyd and L.~Vandenberghe, \emph{Convex Optimization}.\hskip 1em plus 0.5em
  minus 0.4em\relax New York, NY: Cambridge University Press, 2004.

\bibitem{KR:98}
K.~Rose, ``Deterministic annealing for clustering, compression, classification,
  regression, and related optimization problems,'' \emph{Proceedings of the
  IEEE}, vol.~80, no.~11, pp. 2210--2239, 1998.

\bibitem{RGB:95}
R.~G. Bartle, \emph{The Elements of Integration and Lebesgue Measure},
  1st~ed.\hskip 1em plus 0.5em minus 0.4em\relax Wiley-Interscience, 1995.

\bibitem{JPL:86}
J.~P. LaSalle, \emph{The Stability and Control of Discrete Processes}, ser.
  Applied Mathematical Sciences.\hskip 1em plus 0.5em minus 0.4em\relax New
  York, NY: Springer Verlag, 1986, vol.~62.

\end{thebibliography}
